\input amstex
\documentstyle {amsppt}
\UseAMSsymbols \vsize 18cm \widestnumber\key{ZZZZZ}

\catcode`\@=11
\def\displaylinesno #1{\displ@y\halign{
\hbox to\displaywidth{$\@lign\hfil\displaystyle##\hfil$}&
\llap{$##$}\crcr#1\crcr}}
\def\ldisplaylinesno #1{\displ@y\halign{
\hbox to\displaywidth{$\@lign\hfil\displaystyle##\hfil$}&
\kern-\displaywidth\rlap{$##$} \tabskip\displaywidth\crcr#1\crcr}}
\catcode`\@=12

\refstyle{C}

\let \ol=\overline
\let \ul=\underline

\let \wt=\widehat

\font\nr=eufb7 at 10pt

\font\srm=cmr10 at 7.5pt

\font\ens=msbm10

\font\sens=msbm10 at 7.5pt

\font\main=cmsy10 at 10pt

\font \fin=lasy8 at 15.4 pt
\def \X{\mathop{\hbox{\nr X}^{\hbox{\srm nr}}}\nolimits}
\def \CC{\mathop{\hbox{\ens C}}\nolimits}

\def \o{\mathop{\hbox{\main O}}\nolimits}

\def \s{\mathop{\hbox {\rm s}}\nolimits}

\def \GL{\mathop{\hbox{\rm GL}}\nolimits}

\def \Aut{\mathop{\hbox{\rm Aut}}\nolimits}

\def \Ad{\mathop{\hbox{\rm Ad}}\nolimits}
\def \rg{\mathop{\hbox{\rm rg}}\nolimits}

\def \ord{\mathop{\hbox{\rm ord}}\nolimits}

\def \GL{\mathop{\hbox{\rm GL}}\nolimits}
\def \SL{\mathop{\hbox{\rm SL}}\nolimits}
\def \SO{\mathop{\hbox{\rm SO}}\nolimits}
\def \Gal{\mathop{\hbox{\rm Gal}}\nolimits}
\def \sGal{\mathop{\hbox{\srm Gal}}\nolimits}

\def \Int{\mathop{\hbox{\rm Int}}\nolimits}

\topmatter
\title Orbites unipotentes et p\^oles d'ordre \\ maximal de la fonction $\mu $ de
Harish-Chandra \endtitle

\rightheadtext{Orbites unipotentes et p\^oles d'ordre maximal}
\author Volker Heiermann \endauthor

\address Institut f\"ur Mathematik, Humboldt-Universit\"at zu
Berlin, Unter den Linden 6, 10099 Berlin, Allemagne \endaddress

\email heierman\@mathematik.hu-berlin.de \endemail

\abstract In a previous work, we have shown that a representation
of a $p$-adic group obtained by (normalized) parabolic induction
from an irreducible supercuspidal representation $\sigma $ of a
Levi subgroup $M$ contains a subquotient which is square
integrable, if and only if Harish-Chandra's $\mu $-function has a
pole in $\sigma $ of order equal to the parabolic rank of $M$. The
aim of the present article is to interpret this result in terms of
Langlands' functoriality principle. \rm

\null\noindent{R\'ESUM\'E: Dans un travail ant\'erieur, nous
avions montr\'e que l'induite parabolique (normalis\'ee) d'une
repr\'esentation irr\'eductible cuspidale $\sigma $ d'un
sous-groupe de Levi $M$ d'un groupe $p$-adique contient un
sous-quotient de carr\'e int\'egrable, si et seulement si la
fonction $\mu $ de Harish-Chandra a un p\^ole en $\sigma $ d'ordre
\'egal au rang parabolique de $M$. L'objet de cet article est
d'interpr\'eter ce r\'esultat en termes de fonctorialit\'e de
Langlands.}

\endabstract

\endtopmatter
\document
Soit $G$ le groupe des points rationnels d'un groupe r\'eductif
connexe d\'efini sur un corps local non archim\'edien $F$. Soit
$P$ un sous-groupe parabolique de $G$ et $M$ un facteur de Levi de
$P$. Notons $i_P^G$ le foncteur d'induction parabolique
normalis\'ee. Soit $\sigma $ une repr\'esentation irr\'eductible
cuspidale de $M$. Dans \cite{H2}, nous avons montr\'e que la
repr\'esentation $i_P^G\sigma $ contient un sous-quotient de
carr\'e int\'egrable, si et seulement si $\sigma $ est un p\^ole
d'ordre $\rg _{ss}(G)-\rg _{ss}(M)$ de la fonction $\mu $ de
Harish-Chandra (o\`u $\rg _{ss}$ d\'esigne le rang semi-simple
d\'eploy\'e), et que la restriction de $\sigma $ sur le centre de
$G$ est un caract\`ere unitaire. Le but de cet article est de
traduire ce r\'esultat en termes de fonctorialit\'e de Langlands.

Plus pr\'ecis\'ement, supposons que l'on dispose d'une
correspondance de Langlands pour les repr\'esentations cuspidales
des sous-groupes de Levi de $G$, v\'erifiant certaines
propri\'et\'es suppl\'ementaires: notons $W_F$ le groupe de Weil
de $F$ et $^LH$ le $L$-groupe d'un sous-groupe r\'eductif $H$ de
$G$. On suppose donc en particulier que l'on sache associer \`a
toute repr\'esentation cuspidale d'un sous-groupe de Levi $M$ de
$G$ un homomorphisme \it admissible \rm $W_F\times\SL _2(\Bbb
C)\rightarrow\ ^LM$. Cet homomorphisme doit par ailleurs \^etre
\it discret, \rm ce qui signifie entre autres que son image n'est
pas contenue dans un sous-groupe de Levi propre de $^LM$. Cet
homomorphisme sera alors appel\'e le param\`etre de Langlands de
$\sigma $. (Le lecteur trouvera plus de pr\'ecisions en {\bf
1.5}.)

Si la restriction \`a $\SL _2(\Bbb C)$ du param\`etre de Langlands
de $\sigma $ est triviale, nous supposons de plus la propri\'et\'e
suivante: soit $M_{\alpha }$ un sous-groupe de Levi de $G$ qui
contient $M$ et qui est minimal pour cette propri\'et\'e. Notons
$M_{\alpha }^{der}$ le groupe d\'eriv\'e de $M_{\alpha }$. Soit
$\chi $ un caract\`ere non ramifi\'e de $M$. Supposons que la
restriction de $\chi $ \`a $M\cap M_{\alpha }^{der}$ ne soit pas
un caract\`ere unitaire, alors que celle au centre de $M_{\alpha
}$ le soit. Si la repr\'esentation induite $i_{P\cap M_{\alpha
}}^{M_{\alpha }}(\sigma\otimes\chi )$ est r\'eductible, alors il
est bien connu que cette repr\'esentation poss\`ede un unique
sous-quotient de carr\'e int\'egrable. On d\'eduit du param\`etre
de Langlands pour $\sigma $ naturellement un param\`etre de
Langlands pour ce sous-quotient. Notre hypoth\`ese
suppl\'ementaire exige que ce param\`etre de Langlands soit
discret.

Si la restriction \`a $\SL _2(\Bbb C)$ du param\`etre de Langlands
de $\sigma $ n'est pas triviale, la situation est plus
compliqu\'ee. Nous avons en quelque sorte besoin que $\sigma $ et
son param\`etre de Langlands ne se comportent pas seulement comme
une repr\'esentation cuspidale, mais \'egalement comme une
repr\'esentation de carr\'e int\'egrable induite par une
repr\'esentation cuspidale $\sigma _0$ dont le param\`etre de
Langlands est la restriction \`a $W_F$ du param\`etre de Langlands
de $\sigma $ et que $\sigma _0$ v\'erifie l'hypoth\`ese ci-dessus.

En fait, on va d\'eduire cette derni\`ere hypoth\`ese d'une
propri\'et\'e (conjecturale) de la correspondance de Langlands. En
effet, dans la situation ci-dessus, il est suppos\'e qu'il existe
une forme int\'erieure $G'$ de $G$ telle que l'on puisse associer
\`a $\sigma $ par la fonctorialit\'e de Langlands (le $L$-paquet
d') une repr\'esentation de carr\'e int\'egrable non cuspidale
$\sigma '$ d'un certain sous-groupe de Levi de $G'$. Il est alors
attendu que $\sigma '$ soit le sous-quotient d'une
repr\'esentation induite par une repr\'esentation cuspidale
$\sigma _0'$ dont le param\`etre de Langlands est la restriction
\`a $W_F$ du param\`etre de Langlands de $\sigma $ et que $\sigma
_0'$ v\'erifie l'hypoth\`ese ci-dessus. Nous renvoyons le lecteur
\`a {\bf 4} et {\bf 5.1} pour plus de pr\'ecisions.

Inversement, nous \'emettons une hypoth\`ese, quand le $L$-paquet
 associ\'e \`a un homomorphisme admissible discret devrait
\^etre form\'e de repr\'esentations cuspidales qui sera
justifi\'ee post\'erieurement.

Ces hypoth\`eses faites, nous associons \`a une repr\'esentation
de carr\'e int\'egrable $\pi $ de $G$ un homomorphisme admissible
discret $\psi _{\pi }:W_F\times\SL _2(\Bbb C)\rightarrow\ ^LG$ qui
est, \`a \'equivalence pr\`es, uniquement d\'etermin\'e par $\pi
$. Nous montrons que l'on obtient ainsi tous les homomorphismes
admissibles qui sont discrets. A l'aide de la classification de
Langlands, nous en d\'eduisons un proc\'ed\'e qui associe \`a
chaque homomorphisme admissible $W_F\times\SL _2(\Bbb
C)\rightarrow\ ^LG$ (le $L$-paquet d') une repr\'esentation
irr\'eductible lisse de $G$,  ainsi qu'un proc\'ed\'e inverse.

Ce travail est donc une sorte de g\'en\'eralisation \`a un groupe
r\'eductif $p$-adique quelconque d'une portion des r\'esultats
classiques de Bernstein-Zelevinsky \cite{R} pour $\GL _n$. Nous
esp\'erons en cons\'equence que les r\'esultats \'enonc\'es ici
seront utiles dans ce sens dans l'avenir.

Remarquons toutefois que nous n'obtenons aucune information sur
les caract\`eres non ramifi\'es $\chi $ de $M$, tels que l'induite
parabolique de $\sigma\otimes\chi $ soit r\'eductible. Dans les
cas o\`u cette information est disponible, en particulier si $G$
est quasi-d\'esploy\'e et la repr\'esentation $\sigma $ est
g\'en\'erique \cite{S}, il est pensable qu'une partie des
hypoth\`eses \'emises ici peut \^etre v\'erifi\'ee. Certains
r\'esultats \'enonc\'es seraient alors des th\'eor\`emes
inconditionn\'es dans ces cas. Mais cela n'est pas encore \'ecrit.
Dans les cas des repr\'esentations lisses irr\'eductibles
g\'en\'eriques des groupes classiques d\'eploy\'es, on dispose
toutefois des r\'esultats de fonctorialit\'e \cite{CKPS} qui,
ensemble avec la correspondance de Langlands pour $\GL _n$
\cite{HT}, permettent d'associer \`a toute repr\'esentation lisse
irr\'eductible g\'en\'erique d'un tel groupe un param\`etre de
Langlands, l'image de cette correspondance n'\'etant connue dans
le cas local pour l'instant que si $G=\SO_{2n+1}(F)$ \cite{JS}.

Le plan de l'article est le suivant: au paragraphe {\bf 0}, nous
introduisons nos notions de base et rappelons le r\'esultat
principal de \cite{H2}. Au paragraphe {\bf 1}, nous revoyons la
notion de $L$-groupe, ses propri\'et\'es ainsi que les notions
d'homomorphisme admissible et de param\`etre de Langlands. Le
paragraphe {\bf 2} est consacr\'e \`a la th\'eorie des $\SL
_2$-triplets. Au paragraphe {\bf 3}, on r\'esume la th\'eorie des
orbites unipotentes d'un groupe r\'eductif connexe complexe, on
rappelle la notion de $L^2$-paire de G. Lusztig, on introduit la
notion d'\'el\'ement semi-simple $q$-distingu\'e et on prouve
finalement que tout \'el\'ement semi-simple $q$-distingu\'e peut
\^etre compl\'et\'e en une $L^2$-paire. Au paragraphe {\bf 4} nous
\'enon\c cons les hypoth\`eses qui sont n\'ecessaires pour
associer \`a un homomorphisme admissible discret une
repr\'esentation de carr\'e int\'egrable, et nous donnons quelques
propri\'et\'es des homomorphismes admissibles. Les r\'esultats
principaux de l'article se trouvent alors dans le paragraphe {\bf
5}.

L'auteur a s\'ejourn\'e durant des \'etapes pr\'eliminaires \`a ce
travail \`a l'Institut for Advanced Study \`a Princeton
(cofinanc\'e par une bourse Feodor-Lynen de la fondation Humboldt
et la bourse DMS 97-29992 de la NSF) et \`a l'IH\'ES \`a
Bures-sur-Yvette.

Il remercie particuli\`erement G. Lusztig, G. Prasad et F. Shahidi
pour des discussions  sur diff\'erents aspects de l'article,
l'Universit\'e Purdue pour son hospitalit\'e lors de
l'aboutissement de ce travail et A.-M. Aubert pour quelques
remarques suppl\'ementaires.

\null \null {\bf 0. Notations et pr\'eliminaires:} Une bonne
partie des notations et conventions introduites ci-dessous devrait
\^etre standard.

\null {\bf 0.1} L'ensemble des caract\`eres rationnels d'un groupe
alg\'ebrique $\ul{H}$ sera not\'e $X^*(\ul{H})$ et celui des
cocaract\`eres $X_*(\ul{H})$. On \'ecrira $\ul{H}^{\circ }$ pour
la composante neutre, $C(\ul{H})$ pour le centre de $\ul{H}$ et
$C_{\ul{H}}(h)$ pour le centralisateur d'un \'el\'ement $h$ de
$\ul{H}$. Si on parle du centralisateur d'un sous-groupe $H_1$ de
$\ul{H}$, on entendra par l\`a l'intersection des $C_{\ul{H}}(h)$,
$h$ parcourant $H_1$.

Si $\ul{H}$ est un groupe r\'eductif et $\ul{T}$ un tore maximal
de $\ul{H}$, $\Sigma (\ul{T})$ d\'esignera l'ensem-ble des racines
de $\ul{T}$ dans l'alg\`ebre de Lie de $\ul{H}$. On notera $Ad $
l'action adjointe de $\ul{T}$ sur l'alg\`ebre de Lie de $\ul {H}$.

Le tore maximal dans le centre de $\ul {H}$ sera not\'e
$T_{\ul{H}}$.

\null {\bf 0.2} Soit $\Cal {G}$ un groupe alg\'ebrique complexe et
$s$ un \'el\'ement semi-simple de $\Cal {G}$. Fixons un
sous-groupe diagonalisable minimal $D$ de $\Cal {G}$ qui contient
$s$. Notons $D_c$ (resp. $D_v$) le sous-groupe de $D$ form\'e des
\'el\'ements $g$ de $D$ tels que, pour tout caract\`ere
alg\'ebrique $\chi :D\rightarrow \Bbb C$, $\vert\chi (g)\vert =1$
(resp. $\chi (g)\in\Bbb R^{>0}$). Alors $D$ est canoniquement
isomorphe au produit direct de $D_v$ et $D_c$. Notons $s_c$ (resp.
$s_v$) la projection de $s$ sur $D_c$ (resp. $D_v$). On appellera
$s_c$ la partie compacte et $s_v$ la partie hyperbolique de $\Cal
{G}$. Remarquons que la d\'ecomposition \it polaire \rm $s=s_cs_v$
est invariante par homomorphisme de groupes alg\'ebriques.

\null {\bf 0.3} Dans tout ce travail $F$ d\'esignera un corps
local non archim\'edien, $q$ le cardinal de son corps r\'esiduel,
$\vert . \vert _F$ sa valeur absolue normalis\'ee, $v_F$ la
valuation discr\`ete associ\'ee et $W_F$ son groupe de Weil. On
fixera un Frob\'enius g\'eom\'etrique $Fr$ dans $W_F$. En
identifiant l'ab\'elianis\'e de $W_F$ avec le groupe multiplicatif
de $F$ par la th\'eorie du corps de classes (normalis\'ee de sorte
que les uniformisantes correspondent aux automorphismes de
Frob\'enius g\'eom\'etriques), on d\'efinit $v_F(\gamma )$ pour un
\'el\'ement $\gamma $ de $W_F$. En particulier, $v_F(Fr)=1$.

\null {\bf 0.4} Le symbole $G$ d\'esignera le groupe des points
rationnels d'un groupe r\'eductif connexe $\underline {G}$
d\'efini sur $F$. On se fixe un tore $F$-d\'eploy\'e maximal
$\underline {A_0}$ dans $\underline{G}$ et un tore maximal
$\underline {T}$ de $\underline {G}$ d\'efini sur $F$ qui contient
$\underline {A_0}$. Fixons \'egalement un sous-groupe parabolique
minimal $\underline {P_0}$ de $\underline {G}$ d\'efini sur $F$ et
contenant $\underline {T}$. On notera $\ul {M_0}$ l'unique facteur
de Levi de $\ul {P_0}$ d\'efini sur $F$ et qui contient $\ul
{A_0}$.

\null {\bf 0.5} On appellera sous-groupe parabolique semi-standard
de $G$ tout groupe qui est le nombre des points rationnels d'un
sous-groupe parabolique de $\underline {G}$ d\'efini sur $F$ et
qui contient $A_0$. S'il contient en outre $P_0$, on l'appellera
un sous-groupe parabolique standard. Dans les deux cas, il existe
un unique facteur de Levi d\'efini sur $F$ qui contient $A_0$. En
\'ecrivant "$P=MU$ est un sous-groupe parabolique (semi-)standard
de $G$", on sous-entend que $U$ est le radical unipotent de $P$ et
que $M$ contient $M_0$. Un tel sous-groupe $M$ de $G$ sera plus
simplement appel\'e sous-groupe de Levi (semi-) standard.

\null {\bf 0.6} Si $M$ est un sous-groupe de Levi semi-standard de
$G$, on notera $A_M$ le tore d\'eploy\'e maximal dans le centre de
$M$ et $\Sigma (A_M)$ l'ensemble des racines pour l'action de
$A_M$ dans l'alg\`ebre de Lie de $M$. On pose
$a_M^*=X^*(A_M)\otimes_{\Bbb Z}\Bbb R$, et on note $a_M$ l'espace
dual. Si $M'$ est un sous-groupe de Levi semi-standard de $G$ qui
contient $M$, on note $a_M^{M'*}$ le sous-espace de $a_M^*$,
annul\'e par $a_{M'}$. On a donc une d\'ecomposition
$a_M^*=a_{M'}^*\oplus a_M^{M'*}$. On notera $\lambda =\lambda
_{M'}+\lambda _M^{M'}$ la d\'ecomposition d'un \'el\'ement
$\lambda $ de $a_M^*$ selon cette d\'ecomposition.

Le choix d'un sous-groupe parabolique semi-standard $P$ de $G$
dont $M$ est un facteur de Levi est \'equivalent au choix d'un
certain ordre $>_P$ sur $a_M^*$. On notera $\Sigma (P)$ l'ensemble
des racines dans $\Sigma (A_M)$ qui sont positives pour $P$.
Rappelons que l'on a une bijection $\alpha\mapsto M_{\alpha }$
entre l'ensemble des racines r\'eduites dans $\Sigma (P)$ et celui
des sous-groupes de Levi semi-standard minimaux contenant $M$. (Le
sous-groupe de Levi $M_{\alpha }$ est le centralisateur de $(\ker
\alpha)^{\circ }$.)

On a un proc\'ed\'e qui associe \`a un \'el\'ement $\lambda $ du
complexifi\'e $a_{M,\Bbb C}^*$ de $a_M^*$ un certain
(quasi-)caract\`ere $\chi _{\lambda }$ de $M$. Si $\alpha $ est un
caract\`ere $F$-rationnel de $M$, et si $\lambda $ est un nombre
complexe, alors $\chi _{\lambda\alpha }(m):=\vert\alpha (m)\vert
_F^{\lambda }$. Un tel caract\`ere de $M$ sera appel\'e \it
caract\`ere non ramifi\'e, \rm et le groupe form\'e de ces
caract\`eres not\'e $\X (M)$. Ce groupe a la structure d'une
vari\'et\'e alg\'ebrique complexe.

\null {\bf 0.7} Fixons un sous-groupe parabolique semi-standard
$P=MU$ de $G$. Le groupe $\X (M)$ agit sur l'ensemble des (classes
d'\'equivalence) de repr\'esentations cuspidales de $M$. Une
orbite $\o $ pour cette action est une vari\'et\'e alg\'ebrique
complexe. Sur $\o $, on a donc une notion de fonction rationnelle.
La fonction $\mu $ de Harish-Chandra est une certaine fonction
rationnelle $\mu =\mu ^G$ d\'efinie sur une telle orbite $\o $ qui
appara\^\i t naturellement dans la formule de Plancherel d'un
groupe $p$-adique \cite{W}.

Si $M$ est un sous-groupe de Levi maximal de $G$ et $\sigma $ une
repr\'esentation cuspidale unitaire de $M$, alors la fonction $\mu
$ ne d\'epend que d'une seule variable complexe que l'on peut
identifier avec $\sigma\otimes\chi _{\lambda }$, $\lambda\in
a_{M,\Bbb C}^{G*}$. On sait que $\sigma\otimes\chi _{\lambda }$,
$\lambda\in a_M^{G*}$, ne peut \^etre qu'un p\^ole de $\mu $, si
$\mu (\sigma )=0$. Inversement, si $\mu (\sigma )=0$, il existe un
unique $\lambda\in a_M^{G*}$, au signe pr\`es, tel que $\mu $ soit
singulier en $\sigma\otimes\chi _{\lambda }$. Ces p\^oles sont
d'ordre $1$, alors que les z\'eros sont d'ordre $2$. Pour que la
repr\'esentation induite $i_P^G(\sigma\otimes\chi _{\lambda })$
soit r\'eductible avec $\lambda$ non nul dans $a_M^{G*}$, il faut
et il suffit que $\sigma\otimes\chi _{\lambda }$ soit un p\^ole de
$\mu $. On trouvera ces r\'esultats dus \`a Harish-Chandra et en
partie \`a A. Silberger dans \cite{Si0} et \cite{Si1} (comparer
\'egalement la remarque \cite{H2, 4.1} relative \`a la
simplicit\'e des p\^oles de la fonction $\mu $).

Si $M$ n'est plus un sous-groupe de Levi maximal de $G$, on a une
formule du produit $\mu =\prod _{\alpha }\mu ^{M_{\alpha }}$,
$\alpha $ parcourant les racines r\'eduites dans $\Sigma (P)$. Ces
facteurs ne d\'ependent que d'une seule variable complexe, $M$
\'etant un sous-groupe de Levi maximal de $M_{\alpha }$. L'ordre
du p\^ole de $\mu $ en un point $\sigma\otimes\chi $ de $\o $ est
alors (par d\'efinition si on veut) \'egal \`a la somme des
p\^oles des diff\'erents facteurs dans la formule du produit. On
le note $\ord _{\sigma\otimes\chi }\mu $. Si ce nombre est
n\'egatif, on parlera d'un z\'ero.

Rappelons que le rang parabolique de $M$ est d\'efini par $\rg
_{ss}(G)-\rg_{ss}(M)$ et signalons le th\'eor\`eme suivant qui est
le r\'esultat principal de \cite{H2}:

\null {\bf Th\'eor\`eme:} \it Soit $P=MU$ un sous-groupe
parabolique de $G$ et $\sigma $ une repr\'esenta-tion
irr\'eductible cuspidale de $M$. Pour que la repr\'esentation
induite $i_P^G\sigma $ poss\`ede un sous-quotient de carr\'e
int\'egrable, il faut et il suffit que $\sigma $ soit un p\^ole de
la fonction $\mu $ de Harish-Chandra d'ordre \'egal au rang
parabolique de $M$ et que la restriction de $\sigma $ \`a $A_G$
soit un caract\`ere unitaire. \rm

\null \null {\bf 1.} On va r\'esumer ci-dessous quelques
propri\'et\'es du $L$-groupe, et on \'enoncera les conjectures
locales de Langlands. Le lecteur trouvera plus des d\'etails dans
\cite{B}. Seule la section {\bf 1.4} n'est peut-\^etre pas
recouverte par la litt\'erature.

\null {\bf 1.1} Soit $\underline {B}$ un sous-groupe de Borel de
$\underline {G}$ contenu dans $\underline {P_0}$ et dont
$\underline {T}$ est un sous-groupe de Levi. Notons
$X^*=X^*(\underline {T})$ le groupe des caract\`eres rationnels de
$\underline {T}$ et $X_*=X_*(\underline {T})$ celui des
cocaract\`eres. Au choix de $\underline {B}$ correspond une base
$\underline {\Delta }$ de $\Sigma (\underline {T})$. Notons $\Phi
(\ul{G})=(X^*,\ul {\Delta }, X_*, \ul {\Delta }^{\vee })$ la
donn\'ee radicielle basique associ\'ee \`a $\ul {B}$ et $\ul {T}$,
et $\Phi ^{\vee }(\ul {G})=(X_*,\ul {\Delta }^{\vee }, X^*, \ul
{\Delta })$ la donn\'ee radicielle basique duale. Il existe - \`a
isomorphisme pr\`es - un unique triplet $(\wt {G}, \wt {B}, \wt
{T})$ form\'e d'un groupe r\'eductif connexe complexe $\wt {G}$,
d'un sous-groupe de Borel $\wt {B}$ de $\wt {G}$ et d'un tore
maximal $\wt {T}$ de $\wt {B}$ de donn\'ee radicielle basique
$\Phi ^{\vee }(\ul {G})$.

\null {\bf 1.2} Fixons une cl\^oture s\'eparable $F^{sep}$ de $F$.
Si $\gamma $ est un \'el\'ement de $\Gamma =\Gal (F^{sep}/F)$,
alors il existe $g\in \ul {G}(F^{sep})$, tel que $g^{\gamma }\ul
{T}g^{-1}=\ul {T}$ et que $g^{\gamma }\ul {B}g^{-1}=\ul {B}$.
L'\'el\'ement $g$ est d\'etermin\'e \`a un multiple par un
\'el\'ement de $\ul {T}(F^{sep})$ pr\`es. On en d\'eduit un
automorphisme $\nu _{\ul {G}}: \Gamma \rightarrow \Aut (\Phi
^{\vee }(G))$. Remarquons que, si $\ul {G'}$ est un autre groupe
r\'eductif connexe d\'efini sur $F$, alors $\nu _{\ul{G}}=\nu
_{\ul{G'}}$, si et seulement si $\ul {G'}$ est une forme
int\'erieure de $\ul {G}$.

Comme $\Aut (\Phi (\ul {G}))=\Aut (\Phi ^{\vee }(\ul {G}))$, on
d\'eduit de tout monomorphisme $\Aut (\Phi ^{\vee }$ $(\ul
{G}))\rightarrow \Aut (\wt {G}, \wt {B}, \wt {T})$ une action $\nu
_{\ul {G}}$ de $\Gamma $ sur $\wt {G}$ que l'on notera $g\mapsto\
^{\gamma }g$. Cette action se factorise par le groupe de Galois de
l'extension galoisienne minimale $K/F$ sur laquelle $\ul {G}$ se
d\'eploie. Cette extension est de degr\'e fini.  On d\'efinit
alors $^LG$ comme \'etant le produit semi-direct $\wt
{G}\rtimes\Gal (K/F)$ d\'eduit de cette action de $\Gal (K/F)$,
i.e., pour tout $g\in\wt {G}$ et tout $\gamma\in\Gal (K/F)$, on a
$\gamma g=\ ^{\gamma }g \gamma $. Cette d\'efinition d\'epend du
choix du monomorphisme $\Aut (\Phi ^{\vee }(\ul {G}))\rightarrow
\Aut (\wt {G}, \wt {B}, \wt {T})$. Si on change ce monomorphisme,
on ne change $\nu _{\ul {G}}$ toutefois que par un automorphisme
int\'erieur $\Int (t)$, $t\in\ul {T}$. Le groupe $^LG$ est donc
d\'etermin\'e \`a isomorphisme int\'erieur pr\`es. C'est un groupe
r\'eductif \it complexe, \rm qui est connexe si et seulement si
$G$ est d\'eploy\'e. \rm

\null {\bf 1.3} Un sous-groupe parabolique de $^LG$ est un
sous-groupe ferm\'e de $G$ qui est \'egal au normalisateur d'un
sous-groupe parabolique $\widehat {P}$ de $\widehat {G}$ et dont
la projection sur le deuxi\`eme facteur est $\Gal (K/F)$. On dit
qu'il est standard, s'il contient $^LB:=\widehat {B}\rtimes \Gal
(K/F)$. Il est alors de la forme $\widehat {P}\rtimes \Gal (K/F)$
avec $\widehat {P}\supseteq\widehat {B}$ (en particulier $\widehat
{P}$ poss\`ede un sous-groupe de Levi $\widehat
{M}\supseteq\widehat {T}$).

Tout sous-groupe parabolique de $^LG$ est conjugu\'e \`a un unique
sous-groupe parabolique standard. Ainsi les classes de conjugaison
des sous-groupes parabo-liques de $^LG$ sont en bijection avec les
classes de conjugaison des sous-groupes paraboliques de $\ul {G}$
stables pour l'action par $\Gal (K/F)$.

On obtient une injection de l'ensemble des sous-groupes
paraboliques standard de $G$ dans l'ensemble des sous-groupes
paraboliques standard de $^LG$, en associant \`a un sous-groupe
parabolique standard $P$ de $G$ le sous-groupe parabolique
$\widehat {P}\rtimes\Gal (K/F)$ de $^LG$. On notera $^LP$ l'image
de $P$.

Cette injection est bijective si et seulement si $G$ est
quasi-d\'eploy\'e. Un sous-groupe parabolique de $^LG$ est dit \it
admissible, \rm s'il est conjugu\'e \`a un sous-groupe parabolique
de la forme $^LP$ avec $P$ un sous-groupe parabolique standard de
$G$. On fait une d\'efinition analogue pour un sous-groupe de Levi
de $^LG$.

\null {\bf 1.4} Remarqons que le tore maximal $T_{\underline G}$
dans le centre $C(\underline {G})$ de $\underline {G}$ est
d\'efini sur $T$. Son groupe des points $F$-rationnel est donc
\'egal \`a $T_G$. Comme le quotient de $\underline{G}$ par son
groupe d\'eriv\'e est isomorphe au quotient de $C(\underline {G})$
par un sous-groupe fini, on a un morphisme surjectif de groupes
alg\'ebriques $\underline{G}\longrightarrow T_{\underline G}$. Des
propri\'et\'es de fonctorialit\'e des $L$-groupes \cite{B, 2.5},
on d\'eduit une inclusion $^LT_G\hookrightarrow\ ^LG$, o\`u
$^LT_G=T_{\widehat {G}}\rtimes\Gal (K/F)$. La composante neutre du
centre de $^LG$ est donc \'egale \`a la composante neutre du
groupe des points fixes de $T_{\widehat {G}}$ pour l'action de
$\Gal (K/F)$. Sa donn\'ee de racines basique est
$(X^*(T_{\underline G})^{\sGal (K/F)}, \{1\}, X_*(T_{\underline
G})^{\sGal (K/F)}, \{1\})$. C'est le $L$-groupe du tore
$F$-d\'eploy\'e maximal $\underline {A}_G$ dans le centre de
$\underline {G}$. En particulier, $A_G$ et $C(\ ^LG)^{\circ}$ ont
m\^eme rang.

Plus g\'en\'eralement, on en d\'eduit que, si $^LM$ est un
sous-groupe de Levi semi-standard de $^LG$ (i.e. $^LM$ contient
$^LT$), la composante neutre du centre de $^LM$ correspond \`a un
certain sous-tore de $A_0$ dont le centralisateur est un
sous-groupe de Levi semi-standard $M$ de $G$. Ainsi, on obtient
une bijection entre les sous-groupes de Levi semi-standard de
$^LG$ et de $G$. En particulier, si $M$ est un tel sous-groupe de
Levi de $G$ et $\alpha\in\Sigma (A_M)$, alors $^LM_{\alpha }$ est
le centralisateur de la composante neutre du noyau de la racine
$\alpha ^{\vee }$ de $\wt {T}$ (remarquons que la coracine $\alpha
^{\vee }$ de $A_M$ a \'et\'e d\'efinie dans \cite{H2, 1.2}). On a
$^LM_{\alpha }=\wt {M}_{\alpha }\rtimes \Gal (K/F)$.

En particulier, l'ensemble des racines $\Sigma (A_M)$ s'identifie
\`a un ensemble de cora-cines pour $T_{^LM}$.

Nous d\'efinissons le rang semi-simple \it d\'eploy\'e \rm $\rg
_{ss}(^LM)$ d'un sous-groupe de Levi semi-standard $^LM$ de $^LG$
par $\rg _{ss}(^LM)=\rg(T_{^LT})-rg(T_{^LM})$. On appellera la
diff\'erence $\rg _{ss}(^LG)-\rg_{ss}(^LM)$ le \it rang
parabolique de $^LM$. \rm  Le r\'esultat suivant est une
cons\'equence imm\'ediate des remarques ci-dessus:

\null {\bf Proposition:} \it Pour tout sous-groupe de Levi
semi-standard $M$ de $G$, $^LM$ et $M$ ont m\^eme rang
parabolique. \rm

\null {\bf 1.5} Un morphisme de groupes $\psi : W_F\times
\SL_2(\CC )\rightarrow\ ^LG$ est dit \it admissible, \rm si

(i) la composition de $\psi $ avec la projection de $^LG$ sur
$\Gal (K/F)$ est \'egale \`a la projection $W_F\rightarrow\Gal
(K/F)$;

(ii) la restriction de $\psi $ \`a $W_F$ est un morphisme de
groupes continu;

(iii) la restriction de $\psi $ \`a $\SL _2(\CC )$ est un
morphisme de groupes alg\'ebriques;

(iv) pour tout $\gamma\in W_F$, $\psi (\gamma ,1)$ est un
\'el\'ement semi-simple de $^LG$;

(v) si $\psi $ se factorise par un sous-groupe de Levi de $^LG$,
alors celui-ci est admissible.

\null Un morphisme admissible $\psi :W_F\times \SL _2(\CC
)\rightarrow\ ^LG$ est dit \it discret, \rm si $\psi (W_F)$ est
relativement compact et si l'image de $\psi $ n'est contenue dans
aucun sous-groupe de Levi propre de $^LG$. Deux morphismes
admissibles $\psi _1, \psi _2 :W_F\times \SL _2(\CC )\rightarrow\
^LG$ sont dits \it \'equivalents, \rm s'il existe $g\in \
\widehat{G}$ avec $\Int (g)\circ \psi _1=\psi _2$, o\`u $\Int (g)$
d\'esigne l'automorphisme int\'erieur de $\widehat {G}$ associ\'e
\`a $g$.

L'homomorphisme admissible $\psi $ est dit \it non ramifi\'e, \rm
si sa restriction au sous-groupe d'inertie de $W_F$ est trivial.

\null {\bf 1.6 Conjecture} (Conjecture locale de Langlands) \it
Les (classes d'\'equivalence de) repr\'esentations irr\'eductibles
lisses de $G$ sont en bijection avec les classes d'\'equiva-lence
de couples $(\psi ,\phi )$ form\'es d'un homomorphisme admissible
$\psi :W_F\times \SL_2 (\CC )\rightarrow\ ^LG$ et d'une
repr\'esentation irr\'eductible $\phi $ du groupe fini $S_{\psi
}:=C_{^LG}$ $(im(\psi ))$ $/C_{^LG}$ $(im(\psi ))^{\circ }$, les
repr\'esentations de carr\'e int\'egrable de $G$ correspondant aux
couples $(\psi ,\phi )$ avec $\psi $ discret. \rm

\null \it Remarque: \rm Ces conjectures ne sont connues en toute
g\'en\'eralit\'e que pour le groupe $\GL _n$ gr\^ace aux travaux
de M. Harris et R. Taylor \cite{HT} (qui \'etablissent la
correspondance pour les repr\'esentations cuspidales, le
prolongement \`a l'ensemble (des classes d'\'equivalence) des
repr\'esentations irr\'eductibles lisses r\'esultant de travaux
pr\'eliminaires de J. Bernstein et A. Zelevinski \cite{R}). Dans
le cas o\`u $G$ est simple d\'eploy\'e de type adjoint, les
repr\'esentations irr\'eductibles lisses correspondant \`a des
couples $(\psi ,\phi )$ avec $\psi $ non ramifi\'e ont \'et\'e
d\'etermin\'ees par G. Lusztig \cite{L2}.

Par ailleurs, comme d\'ej\`a signal\'e dans l'introduction, la
correspondance de Langlands pour les repr\'esentations
irr\'eductibles lisses g\'en\'eriques des groupes classiques
d\'eploy\'es est une cons\'equence des r\'esultats de
fonctorialit\'e dus \`a J.W. Cogdell, H.H. Kim, I.I.
Piatetski-Shapiro et F. Shahidi \cite{CKPS}, joints \`a la
correspondance de Langlands pour $\GL _n$, l'image de cette
correspondance n'\'etant toutefois connue pour l'instant que si
$G=\SO_{2n+1}(F)$ \cite{JS}.

\null {\bf 1.7} \it D\'efinition: \rm Soit $\pi $ une
repr\'esentation irr\'eductible lisse de $G$. Si $\pi $ correspond
au couple $(\psi ,\phi )$ par les conjectures locales de
Langlands, alors on appellera $\psi $ le \it param\`etre de
Langlands \rm de $\pi $. Deux repr\'esentations irr\'eductibles
lisses de $G$ sont dans le m\^eme \it L-paquet \rm si et seulement
s'ils ont m\^eme param\`etre de Langlands.

\null {\bf 1.8} Dans le cas des tores, la correspondance
conjectur\'ee en {\bf 1.6}, a \'et\'e prouv\'ee par Langlands dans
\cite {L}. Si $G$ est un tore d\'eploy\'e, $G=(F^{\times })^d$,
elle associe au caract\`ere non ramifi\'e $\chi :G\rightarrow \Bbb
C$, $(x_1,\cdots ,x_d)\mapsto \vert x_1\vert _F^{\lambda _1}
\cdots \vert x_d\vert _F^{\lambda _d}$, $\lambda _1,\cdots
,\lambda _d\in \Bbb C$, l'homomorphisme non ramifi\'e
$W_F\rightarrow\ ^LG=(\Bbb C^{\times })^d$ qui envoie $Fr$ sur
$s:=(q^{\lambda _1},\dots ,q^{\lambda _d})$. On dira que $s$
correspond \`a $\chi $ (ou \`a $\lambda $, si $\chi =\chi
_{\lambda }$ pour un $\lambda\in a_{G,\Bbb C}^*$) par la \it
correspondance de Langlands pour les tores, \rm et vice-versa. Si
$\lambda $ est r\'eel (i.e. $\lambda\in a_G^*$), alors la partie
compacte dans la d\'ecomposition polaire de $s$ (cf. {\bf 0.2})
est triviale.

D'autre part, si $M$ est un sous-groupe de Levi semi-standard d'un
groupe $p$-adique et si $\alpha $ est une racine relative \`a
$A_M$, alors le caract\`ere non ramifi\'e $\chi _{\alpha }$ de $M$
correspond \`a l'\'el\'ement semi-simple $\alpha (q)$ de $T_{^LM}$
(o\`u on consid\`ere $\alpha $ comme coracine relative \`a
$T_{^LM}$ (cf. {\bf 1.4})).

\null {\bf 2.} Dans cette section, on se fixe un groupe
alg\'ebrique connexe complexe semi-simple $\Cal G$ d'alg\`ebre de
Lie $\frak g$. Les r\'esultats r\'esum\'es dans {\bf 2.1}-{\bf
2.3} sont bien connus, et on renvoie le lecteur \`a [Ca] pour plus
de pr\'ecisions.

\null {\bf 2.1 Th\'eor\`eme:} (Jacobson-Morozow) \it Soit $N$ un
\'el\'ement nilpotent non nul de $\frak g$. Alors il existe un
morphisme d'alg\`ebres de Lie $\phi : sl _2\rightarrow \frak g$
v\'erifiant $\phi (\pmatrix 0 & 1 \\ 0 & 0 \endpmatrix)=N$.

\null {\bf 2.2} \it D\'efinition: \rm Un triplet $\{ H, N, \ol
{N}\}$ dans $\frak g$ tel qu'il existe un morphisme d'alg\`ebres
de Lie $\phi: \frak {sl} _2\rightarrow\frak g$ v\'erifiant $\phi
(\pmatrix 0 & 1 \\ 0 & 0 \endpmatrix)=N$, $\phi (\pmatrix 1 & 0 \\
0 & -1\endpmatrix)=H$ et $\phi (\pmatrix 1 & 0 \\ 1 & 1
\endpmatrix)=\ol {N}$ est appel\'e un $\frak sl _2$-triplet.

\null {\bf 2.3 Proposition:} \it Une condition n\'ecessaire et
suffisante pour que $\{ H, N, \ol {N}\}$ soit un $\frak sl
_2$-triplet dans $\frak g$, est que ses \'el\'ements v\'erifient
les relations $[H,N]=2N$, $[H,\ol {N}]=-2N$ et $[N,\ol {N}]=H$.

Deux $\frak sl _2$-triplet sont \'egaux, s'ils ont deux
\'el\'ements en commun.\rm

\null En utilisant l'application exponentielle, on d\'eduit alors
de {\bf 2.1} et {\bf 2.3} le corollaire suivant:

\null {\bf 2.4 Corollaire:} \it Soit $s$ un \'el\'ement
semi-simple et $u=\exp(N)$ un \'el\'ement unipotent de $\Cal {G}$.
Supposons que $H={2\over\log q}\log s$ et $N=\log u$ font partie
d'un $\frak sl _2$-triplet dans l'alg\`ebre de Lie de $G$.

Alors il existe un unique morphisme de groupes alg\'ebriques $\SL
_2(\Bbb C) \rightarrow\Cal {G}$ qui envoie $\pmatrix q^{1/2} & 0
\\ 0 & q^{-1/2}\endpmatrix$ sur $s$ et $\pmatrix 1 & 1 \\ 0 &
1\endpmatrix$ sur $\exp (N)$. Son image est un groupe r\'eductif
semi-simple de rang $1$ et l'image de $\pmatrix 1 & 0 \\ 1 & 1
\endpmatrix$ est \'egale au conjugu\'e de $\exp (N)$ par l'unique
\'el\'ement du groupe de Weyl de ce groupe. \rm

\null {\bf 2.5} Une condition n\'ecessaire pour que les
hypoth\`eses du corollaire {\bf 2.4} soient v\'erifi\'ees
relativement \`a un \'el\'ement semi-simple $s$ et un \'el\'ement
unipotent $u$ de $\Cal {G}$ est \'evidemment $sus^{-1}=u^q$. Mais,
cette condition n'est pas suffisante. Cependant, le r\'esultat
suivant vaut:

\null {\bf Proposition:} \it Soient $s$ un \'el\'ement semi-simple
et $u$ un \'el\'ement unipotent de $\Cal {G}$ tels que
$sus^{-1}=u^q$. Alors il existe un \'el\'ement semi-simple $s_1$
dans $\Cal {G}$ et un morphisme de groupes alg\'ebriques $\phi
:\SL _2(\Bbb C)\rightarrow\Cal {G}$ v\'erifiant
$$\phi (\pmatrix q^{1/2} & 0\\ 0 &
q^{-1/2}\endpmatrix)=s_1\qquad\hbox{\rm et}\qquad\phi (\pmatrix 1 & 1 \\
0 & 1\endpmatrix)=u,$$ tels que $ss_1^{-1}$ commute aux
\'el\'ements dans l'image de $\phi $. En particulier, $s$ et $s_1$
commutent.

\null Preuve: \rm  Ce r\'esultat est contenu dans \cite{KL,
paragraph 2}. Cependant ces auteurs font globalement l'hypoth\`ese
que le groupe d\'eriv\'e de $\Cal {G}$ est simplement connexe. On
va reprendre les arguments de \cite{KL} pour montrer que cette
hypoth\`ese n'intervient pas dans ce r\'esultat: choisissons un
homomorphisme de groupes alg\'ebriques $\phi :\SL _2(\Bbb
C)\rightarrow \Cal {G}$ v\'erifiant $\phi (\pmatrix 1 & 1
\\ 0 & 1\endpmatrix )=u$. Posons $\alpha ^{\vee }(z)=\pmatrix z^{1/2} & 0\\ 0 &
z^{-1/2}\endpmatrix$,
$$M(u)=\{(g,z)\in\Cal {G}\times\Bbb {C}^*\vert\
gug^{-1}=u^z\}\qquad\hbox{\rm et}$$ $$M_{\phi }=\{(g,z)\in\Cal
{G}\times\Bbb {C}^*\vert\ g\phi (h)g^{-1}=\phi (\alpha ^{\vee }(z)
h\alpha ^{\vee }(z^{-1}))\ \hbox{\rm pour tout $h\in\SL _2(\Bbb
C)$}\}.$$ Comme $\phi $ est uniquement d\'etermin\'e \`a
conjugaison par un \'el\'ement de $C_{\Cal {G}}(u)$ pr\`es, on
montre comme dans \cite{BV, 2.1} que $M_{\phi }$ est un groupe
r\'eductif maximal de $M(u)$. (Ni la preuve dans \cite{BV}, ni
cette g\'en\'eralisation n'utilisent d'hypoth\`ese sur le groupe
d\'eriv\'e de $\Cal {G}$.)

Rappelons que tout groupe alg\'ebrique est le produit semi-direct
d'un sous-groupe r\'eductif maximal par son radical unipotent et
que deux sous-groupes r\'eductifs maximaux sont conjugu\'es par un
\'el\'ement du radical unipotent. Les sous-groupes r\'eductifs
maximaux de $M(u)$ sont donc n\'ecessairement de la forme $M_{\phi
}$ avec $\phi $ comme ci-dessus. Comme $(s,q)$ est un \'el\'ement
semi-simple de $M(u)$, il existe $\phi :\SL _2(\Bbb C)\rightarrow
\Cal {G}$ v\'erifiant $\phi (\pmatrix 1 & 1\\ 0 & 1\endpmatrix
)=u$, tel que $(s,q)\in M_{\phi }$. Posons $s_1=\phi (\pmatrix
q^{1/2} & 0\\ 0 & q^{-1/2}\endpmatrix )$. Comme $(s_1,q)\in
M_{\phi }$, il en r\'esulte que $(ss_1^{-1},1)\in M_{\phi }$. Par
suite, $ss_1^{-1}$ commute avec les \'el\'ements dans l'image de
$\phi $ et en particulier avec $s_1$.

\hfill {\fin 2}

\null {\bf 2.6 D\'efinition:} Soit $s$ un \'el\'ement semi-simple
de $\Cal {G}$ et $N$ un \'el\'ement nilpotent de $\frak g$. On dit
que $(s,N)$ est une \it $L^2$-paire \rm relative \`a $\Cal {G}$,
si $Ad (s)N=qN$, et si tout tore de $\Cal G$ qui centralise
simultan\'ement $s$ et $N$ est inclus dans le centre de $\Cal
{G}$.

\null \it Remarque: \rm Cette notion est celle introduite par G.
Lusztig dans \cite {L1} pour $\Cal{G}$ un groupe semi-simple, \`a
la diff\'erence pr\`es que Lusztig suppose $Ad (s)N=q^{-1}N$.

\null \null {\bf 3.} Nous continuons \`a noter $\Cal G$ un groupe
r\'eductif complex connexe et $\frak {g}$ son alg\`ebre de Lie.
Les notions et propri\'et\'es donn\'ees dans {\bf 3.1} - {\bf 3.4}
ci-dessous relatives aux orbites nilpotentes et unipotentes d'un
groupe r\'eductif complexe sont bien connues. Le lecteur pourra
consulter par exemple \cite {Ca} pour des preuves d\'etaill\'ees.

\null {\bf 3.1} \it D\'efinition: \rm Un \'el\'ement nilpotent $N$
de $\frak g$ est dit distingu\'e, si et seulement si $\exp (N)$
n'est contenu dans aucun sous-groupe de Levi propre de $\Cal G $.

\null {\bf 3.2 Proposition:} \it Supposons $\Cal G$ semi-simple de
type adjoint. Soit $N$ un \'el\'ement nilpotent non nul dans
$\frak g$ et $\{H, N, \ol {N}\}$ un $\frak sl _2$-triplet
contenant $N$. Posons $\frak g_H (\lambda )=\{Z\in\frak g\vert\
\roman {ad} (H)Z=\lambda Z\}$.

Alors $\dim \frak g _H(0)\geq\dim\frak g_H(2)$. Pour que $N$ soit
distingu\'e dans $G$, il faut et il suffit que $\dim\frak g_H(0)$
$=\dim\frak g_H(2)$. On a alors $\frak g_H(\lambda )=0$ pour tout
entier impair $\lambda $. \rm

\null \it Remarque: \rm Si la derni\`ere propri\'et\'e de la
proposition est v\'erifi\'ee, on dit que $H$ et $\exp (H)$ sont
\it pairs. \rm

\null {\bf 3.3}  Supposons $\Cal G$ simple de type adjoint. Fixons
un sous-groupe de Borel $\Cal{B}=\Cal{T}\Cal{U}$ de $\Cal G$ et
notons $\Delta _1$ la base correspondante du syst\`eme de racines
$\Sigma _1=\Sigma (\Cal{T})$ de $\Cal G$. Soit $J$ un
sous-ensemble de $\Delta _1$ et $n_J: \Sigma _1\rightarrow \Bbb N$
d\'efini par $n_J(\sum_{\alpha \in\Delta _1}\lambda _{\alpha
}\alpha )=2\sum _{\alpha\in J} \lambda_{\alpha }$.

Le sous-groupe parabolique $\Cal P_J$ de $\Cal G $ associ\'e \`a
$J$ est dit \it distingu\'e, \rm si et seulement si
$$\vert\{\beta\in\Sigma _1 \vert n_J(\beta )=2\}\vert = \rg(\Cal
G)+\vert\{\beta \in\Sigma _1\vert n_J(\beta )=0\}\vert .$$

\null {\bf Proposition:} \it Soit $\frak g _J(2)$ (resp. $\frak g
_J(0)$) la somme directe des sous-espaces de $\frak g$ de poids
$\alpha $ v\'erifiant $n_J(\alpha )=2$ (resp. $n_J(\alpha )=0$ et
de l'alg\`ebre de Lie de $\Cal {T}$).

Alors $\dim \frak g_J(0)\geq\dim\frak g_J(2).$

\null {\bf 3.4 Th\'eor\`eme:} (Bala-Carter) \it Supposons $\Cal G$
simple de type adjoint. Soit $\Cal P =\Cal M\Cal U$ un sous-groupe
parabolique distingu\'e de $\Cal G$.

Alors il existe dans l'alg\`ebre de Lie de $\Cal U$ une unique
orbite pour l'action de $\Cal P$ qui est dense. Elle est form\'ee
d'\'el\'ements nilpotents distingu\'es. L'application ainsi
d\'efinie induit une bijection entre les classes de conjugaison de
sous-groupes paraboliques distingu\'es de $\Cal G$ et les classes
de conjugaison d'\'el\'ements nilpotents distingu\'es dans $\frak
g$. \rm

\null {\bf 3.5 D\'efinition:} Soit $s$ un \'el\'ement semi-simple
de $\Cal G$. Notons $\frak g^{der}$ l'alg\`ebre de Lie du groupe
d\'eriv\'e de $\Cal G$. Si $\lambda $ est une valeur propre de $Ad
(s)_{\vert \frak g^{der}}$, notons $\frak g_s(\lambda )$ l'espace
propre associ\'e \`a cette valeur propre. Alors $s$ est dit \it
$q$-distingu\'e, \rm si $\dim(\frak g_s(q))\geq \dim (\frak
g_s(1))$.

\null\it Remarque: \rm Il r\'esultera de la preuve du th\'eor\`eme
{\bf 3.6} (cf. 3.6.2) que l'in\'egalit\'e ci-dessus est en fait
une \'egalit\'e.

\null {\bf 3.6} La preuve du th\'eor\`eme suivant \'etait l'objet
de \cite {H1}. Comme ce manuscrit n'a pas \'et\'e publi\'e, nous
la reproduisons ci-dessous.

\null {\bf Th\'eor\`eme:} \it Soit $s$ un \'el\'ement semi-simple
de $\Cal G$. Pour que $s$ soit $q$-distingu\'e, il faut et il
suffit que $s$ fasse partie d'une $L^2$-paire $(s,N)$ pour $\Cal
G$. La composante nilpotente $N$ est d\'etermin\'ee par $s$ \`a
multiplication par un scalaire non nul pr\`es.

\null Preuve: \rm La preuve r\'esultera des propositions
3.6.1-3.6.3 ci-dessous.

\null Notons $s_vs_c$ la d\'ecomposition polaire de $s$ (cf. {\bf
0.2}).

 \null (3.6.1) \it Soit $s$ un \'el\'ement
$q$-distingu\'e de $\Cal{G}$. Supposons $s=s_v$. Fixons un tore
maximal $\Cal{T}$ de $\Cal{G}$ contenant $s$ et un ensemble
$\Delta _1$ de racines simples relatives \`a $\Cal{T}$ pour lequel
$s$ est positif. Alors toute racine $\alpha $ de $\Delta _1$
v\'erifie ou $\alpha (s)=1$ ou $\alpha (s)=q$. L'in\'egalit\'e
dans {\bf 3.5} est en fait une \'egalit\'e. Il existe un
\'el\'ement nilpotent distingu\'e $N$ dans l'alg\`ebre de Lie de
$\Cal{G}$, tel que $\Ad(s) N=qN$.

\null Preuve de (3.6.1): \rm Notons $\Cal{B}=\Cal{T}\Cal{U}$ le
sous-groupe de Borel de $\Cal{G}$ correspondant \`a $\Delta _1$.
Posons $\s =s^{2\pi i/\log q}$. On a $\s =\s _c$. Le groupe
$\Cal{G'}=C_{\Cal G}(\s_c)^{\circ }$ est r\'eductif connexe de
syst\`eme de racines $\Sigma _1'=\{\alpha \in \Sigma _1\vert \
\alpha (s)=1\}=\{\alpha \in \Sigma _1\vert \ \alpha (s)\in
q^{\hbox {\sens Z}}\}$ relatif \`a $\Cal{T}$. Un sous-groupe de
Borel est donn\'e par $\Cal{B'}=\Cal{B}\cap\Cal{G'}$. Notons
$\Delta _1'$ la base de $\Sigma _1'$ correspondant \`a $\Cal{B'}$.
Posons $J_1=\{ \alpha \in \Delta _1'\vert \alpha (s)=1\}$ et
$J_q=\{ \alpha \in \Delta _1'\vert \alpha (s)=q\}$. On montrera
d'abord que $\Delta _1'=J_1\cup J_q$.

Soit $\Cal{M'}_{\Delta _1'-(J_1\cup J_q)}$ le sous-groupe de
L\'evi standard de $\Cal{G'}$ dont les racines simples sont les
\'el\'ements de $J_1\cup J_q$. On a $s\in\Cal{M'}_{\Delta _1 '-
(J_1\cup J_q)}$ et toute racine positive $\alpha \in \Sigma _1'$
qui v\'erifie $\alpha (s)=q$ ou $\alpha (s)=1$ est combinaison
lin\'eaire d'\'el\'ements de $J_1\cup J_q$.

Le quotient de $\Cal{M'}_{\Delta _1'-(J_1\cup J_q)}$ par son
centre est homomorphe par une bijection \`a un groupe semi-simple
de type adjoint $\Cal{M'}^{ad}$ qui lui est produit direct de
groupes simples de type adjoint, $\Cal{M'}^{ad}=\Cal{M}_1' \times
\cdots\times \Cal{M}_r'$. On a $\vert J_1\cup J_q \vert = \rg
_{ss}(\Cal{M'}_{\Delta '-(J_1\cup J_q)})=\sum _i \rg
_{ss}(\Cal{M}_i')$. Notons $J_{1,i}$ (resp. $J_{q,i}$) le
sous-ensemble de $J_1$ (resp. $J_q$), form\'e de racines pour
$\Cal{M}_i'$, $\Sigma _i'$ la composante irr\'eductible de $\Sigma
'$ correspondant \`a $\Cal{M}_i'$ et $\Cal{P}_{i,J_{1,i}}'$ le
sous-groupe parabolique standard de $\Cal{M}_i'$ de L\'evi
$\Cal{M}_{i,J_{1,i}}'$. On d\'eduit de {\bf 3.3} que $\dim {\frak
g}_s(q)\leq \dim {\frak g}_s(1)$, d'o\`u
$$\eqalign {
\rg _{ss}\Cal{G}&\leq \vert \{\alpha \in \Sigma ^+\vert \alpha
(s)=q\}\vert -
              2 \vert \{\alpha \in \Sigma ^+\vert \alpha (s)=1\}\vert \cr
&\leq \sum _i \vert J_{1,i}\cup J_{q,i}\vert = \vert J_1\cup J_q
\vert \leq \rg _{ss}\Cal{G}. \cr }\leqno \hbox {(*)}$$ On a donc
l'\'egalit\'e partout. En particulier, le syst\`eme de racines
$\Sigma _1'$ a le m\^eme rang que $\Sigma _1$, $J_0\cup J_1=\Delta
_1'$, et l'in\'egalit\'e (*) (et par suite l'in\'egalit\'e dans
{\bf 3.5}) est en fait une \'egalit\'e.

Ceci prouve que le parabolique $\Cal{P}_{J_1}'$ est distingu\'e
dans $\Cal{G}$. Par le th\'eor\`eme de Bala-Carter {\bf 3.4}, on
peut donc trouver un \'el\'ement nilpotent distingu\'e dans
l'alg\`ebre de Lie du radical unipotent de $\Cal{P}_{J_1}'$, tel
que, pour un certain $sl_2$-triplet pour $N$, $s$ soit un
\'el\'ement du tore de rang $1$ de $\Cal{G'}$ d\'eduit de ce
triplet. Comme $\Cal{G}$ et $\Cal{G'}$ ont m\^eme rang et
qu'aucune racine $\alpha $ de $\Sigma\backslash\Sigma '$ ne peut
v\'erifier $\alpha (s)=1$, $N$ est distingu\'e pour $\Cal{G}$ par
{\bf 3.2}. Il est donc pair. Le r\'esultat sur $s$ en suit. \hfill
{\fin 2}

\null (3.6.2) \it On garde les notations et hypoth\`eses du lemme
pr\'ec\'edent, mais on ne suppose plus $s=s_v$.

L'in\'egalit\'e dans {\bf 3.5} est en fait une \'egalit\'e. Il
existe un \'el\'ement nilpotent $N$ dans l'alg\`ebre de Lie de
$\Cal{G}$, tel que $(s,N)$ soit une $L^2$-paire.

\null Preuve de (3.6.2): \rm Soient $\Cal{B}$ et $\Sigma _1$ comme
dans 3.6.1. La composante neutre $\Cal{G'}$ du centralisateur de
$s_c$ dans $\Cal{G}$ est un groupe r\'eductif connexe de syst\`eme
de racines $\Sigma _1'=\{\alpha \in \Sigma \vert \ \alpha
(s_c)=1\}$. Le groupe $\Cal{B'}=\Cal{B}\cap\Cal{G'}$ est un
sous-groupe de Borel de $\Cal{G}$ pour lequel $s$ est positif et
distingu\'e. On a $\alpha (s)=\alpha (s_v)$ pour tout $\alpha \in
\Sigma _1'$. On d\'eduit de 3.6.1 et de son in\'egalit\'e (*) que
$\Cal{G}$ et $\Cal{G'}$ ont m\^eme rang semi-simple et que
l'in\'egalit\'e dans {\bf 3.5} est en fait une \'egalit\'e. En
particulier, le groupe $\Cal{G'}/Z(\Cal{G})$ est semi-simple.

Gr\^ace \`a 3.6.1, on peut choisir un \'el\'ement nilpotent
distingu\'e dans l'alg\`ebre de Lie de $\Cal{G'}$ qui v\'erifie
$\Ad (s)N=qN$. Il reste \`a voir qu'aucun sous-tore non trivial de
$G/Z(G)$ ne peut centraliser simultan\'ement $s$ et $N$. En effet,
comme la d\'ecomposition polaire $s=s_vs_c$ est invariant par
homomorphismes de groupes alg\'ebriques, un tel tore $S$
centraliserait $s_c$. Par suite, $S\subset \Cal{G'}/C(\Cal{G})$.
Comme $S$ centralise \'egalement $N$ et que $N$ est distingu\'e,
on en d\'eduit $S=1$. \hfill {\fin 2}

\null (3.6.3) \it Soit $(s,N)$ une $L^2$-paire relative \`a $q$.
Alors $s$ est $q$-distingu\'e. L'\'el\'ement nilpotent $N$ est
d\'etermin\'e par $s$ \`a multiplication par un scalaire non nul
pr\`es.

\null Preuve: \rm Quitte \`a remplacer $\Cal{G}$ par son quotient
par le centre, on peut supposer $\Cal{G}$ semi-simple. Comme
$(s,N)$ est une $L^2$-paire, le centralisateur connexe $\Cal{G'}$
de la partie compacte de $s$ est semi-simple. En particulier,
$\Cal{G}$ est de m\^eme rang que $\Cal{G'}$. La relation $\Ad
(s)N=qN$ prouve que $N$ est dans l'alg\`ebre de Lie de $\Cal{G'}$.
Comme $(s,N)$ est une $L^2$-paire, $N$ est un \'el\'ement
nilpotent distingu\'e de l'alg\`ebre de Lie de $\Cal{G'}$. Il est
prouv\'e dans \cite{L1, exemple 2.4} que $s_v$ est un \'el\'ement
du tore $S$ de rang $1$ d\'eduit d'un certain $sl_2$-triplet pour
$N$ dans l'alg\`ebre de Lie de $\Cal{G'}$. Il en r\'esulte que $s$
est $q$-distingu\'e.

Choisissons un tore maximal $\Cal{T}$ de $\Cal{G'}$ contenant $S$
et un ensemble $\Delta _1$ de racines simples relatives \`a
$\Cal{T}$ pour lequel $s$ est positif. On d\'eduit de {\bf 3.2}
que toute racine $\alpha $ dans $\Delta _1$ v\'erifie $\alpha
(s)\in \{1,q\}$ et que le parabolique standard $\Cal{P}_J'$,
$J=\{\alpha \vert \alpha (s)=1\}$, est celui associ\'e \`a $N$ par
le th\'eor\`eme de Bala-Carter. Il est donc distingu\'e.
L'unicit\'e de $N$ \`a multiplication par un scalaire pr\`es
r\'esulte alors de la th\'eorie des orbites unipotentes (cf.
\cite{Ca, proposition 5.8.5}). \hfill {\fin 2}

\null {\bf 3.7} En consid\'erant $\Cal {G}$ comme le dual d'un
certain groupe r\'eductif connexe d\'eploy\'e d\'efini sur $F$, on
obtient le r\'esultat suivant:

\null {\bf Proposition:} \it Soit $(s,N)$ une $L^2$-paire relative
\`a $\Cal {G}$. Supposons que $s$ soit un \'el\'ement du groupe
d\'eriv\'e de $\Cal {G}$ et que la partie compacte de $s$ dans sa
d\'ecomposition polaire soit triviale. Alors il existe un
morphisme de groupes alg\'ebriques $\phi :\SL _2(\Bbb
C)$ $\rightarrow \Cal {G}$ qui v\'erifie $$\phi (\pmatrix q^{1/2} & 0\\
0 & q^{-1/2}\endpmatrix)=s\qquad\hbox{\rm et}\qquad\phi(\pmatrix 1 & 1 \\
0 & 1\endpmatrix)=\exp (N).$$ Il est admissible et discret et, \`a
\'equivalence pr\`es, uniquement d\'etermin\'e par $s$.

\null Preuve: \rm Par {\bf 2.5}, on peut choisir $\phi :\SL
_2(\Bbb C)\rightarrow\Cal G$ et un \'el\'ement semi-simple $s_1$
de $\Cal {G}$ tel que $\phi(\pmatrix 1 & 1 \\0 &
1\endpmatrix)=\exp (N)$, $\phi (\pmatrix q^{1/2} & 0\\ 0 &
q^{-1/2}\endpmatrix)=s_1$ et que $z_1:=ss_1^{-1}$ commute aux
\'el\'ements dans l'image de $\phi $. L'\'el\'ement $s_1$
appartient n\'ecessairement au groupe d\'eriv\'e de $\Cal {G}$ et
la partie compacte dans la d\'ecomposition polaire de $s_1$ est
triviale, celle-ci \'etant invariante par morphisme de groupes
alg\'ebriques. Il en est donc de m\^eme de $z_1$. Il suffit de
montrer qu'en fait $z_1=1$ sous nos hypoth\`eses. \'Ecrivons
$N=\sum _{\alpha }N_{\alpha }$, o\`u $N_{\alpha }$ est un
\'el\'ement de l'alg\`ebre de Lie de $\Cal {G}$ de poids $\alpha
$. Notons $n$ le rang semi-simple de $\Cal {G}$. Comme $N$ est
distingu\'e, il doit y \^etre $n$ racines lin\'eairement
ind\'ependantes $\alpha $ telles que $N_{\alpha }\ne 0$. Ces
racines v\'erifient $\alpha (z_1)=1$. Comme $z_1$ est un
\'el\'ement du groupe d\'eriv\'e de $\Cal {G}$ dont la partie
compacte est trivial, ceci prouve que $z_1=1$.

Le morphisme est discret, puisque $(s,N)$ est une $L^2$-paire, et
il est alors imm\'ediat qu'il est admissible. L'unicit\'e vient de
{\bf 2.4}. \hfill{\fin 2}

\null {\bf 4.} On va maintenant introduire les hypoth\`eses dont
nous aurons besoin. Rappelons que le symbole $G$ d\'esigne le
groupe des points $F$-rationnels d'un groupe r\'eductif connexe
d\'efini sur $F$.

\null {\bf 4.1} Le r\'esultat suivant est bien connu. Par manque
de r\'ef\'erence, nous en donnons une preuve:

\null {\bf Lemme:} \it Soit $\psi :W_F\rightarrow\ ^LG$ un
homomorphisme admissible. Alors on peut \'ecrire $\psi =\psi
_{nr}\psi _f$ o\`u $\psi _f$ est un homomorphisme admissible \`a
image finie et $\psi _{nr}$  un homomorphisme admissible non
ramifi\'e (donc trivial sur le groupe d'inertie) \`a valeurs dans
le centralisateur de l'image de $\psi $.

\null Remarque: \rm Le morphisme $\psi _{nr}$ devrait \^etre
li\'e, par les conjectures locales de Langlands {\bf 1.6}, au
caract\`ere central d'une repr\'esentation de $G$ de param\`etre
de Langlands $\psi $.

\null \it Preuve: \rm Le groupe $W_F$ est un produit semi-direct
du sous-groupe d'inertie $I_F$ avec le sous-groupe cyclique
engendr\'e par l'automorphisme de Frob\'enius $Fr$,
$W_F=I_F\rtimes\langle Fr\rangle $. Le groupe $I_F$ \'etant
compact, $\psi (I_F)$ est fini. (C'est une propri\'et\'e bien
connue des topologies profinies et complexes.) Il existe donc une
extension de degr\'e fini $K/F$ telle que $\psi $ se factorise par
$I_K\rtimes \langle Fr\rangle $. Consid\'erons l'homomorphisme
induit $\ol{\psi }: I_F/I_K\rtimes\langle Fr\rangle\rightarrow\
^LG$. L'automorphisme int\'erieur du groupe fini $\ol{\psi
}(I_F/I_K)$ induit par $\ol {\psi }(Fr)$ est \'evidemment d'ordre
fini. Il existe donc un entier $k\geq 1$, tel que $s=\ol {\psi
}(Fr^k)$ centralise $\ol{\psi }(I_F/I_K)$. Par suite, il
centralise \'egalement l'image de $\psi $. Comme le nombre de
composantes connexes d'un groupe alg\'ebrique est fini, il existe
un entier $l$ tel que $s^l\in\widehat {G}$. Le centralisateur
connexe d'un \'el\'ement semi-simple dans un groupe r\'eductif
connexe \'etant r\'eductif, le centralisateur de $s^l$ dans $^LG$
est un groupe r\'eductif $^LG_1$ qui contient l'image de $\psi $.
Notons $D$ un sous-groupe alg\'ebrique diagonalisable contenu dans
le centre de $^LG_1$ et contenant $s^l$. On a $D=D_d\times D_f$,
o\`u $D_f$ est un sous-groupe fini et $D_d$ un tore. Notons $m$
l'ordre de $D_f$. Alors $s^{lm}$ appartient \`a $D_d$. Il existe
donc $s_1\in D_d$ tel que $s_1^{klm}=s^{lm}$. L'homomorphisme
$\psi _f:I_F\rtimes\langle Fr\rangle \rightarrow\ ^LG$,
$(h,Fr^j)\mapsto\psi (h,Fr^j)s_1^{-j}$ est bien d\'efini, continu
et \`a image finie. On en d\'eduit le lemme. \hfill {\fin 2}

\null Le r\'esultat suivant est crucial pour la suite. Il a
d\'ej\`a \'et\'e remarqu\'e dans \cite {K}. Nous pr\'esentons
ci-dessous une preuve l\'eg\`erement diff\'erente.

\null {\bf Proposition:} \it Soit $\psi : W_F\times \SL _2(\Bbb
C)\rightarrow\ ^LG$ un homomorphisme admissible. Alors le
centralisateur de $\psi (W_F)$ dans $^LG$ est un groupe
r\'eductif.

\null Preuve: \rm Par le lemme, on peut \'ecrire $\psi =\psi
_{nr}\psi _f$. Rappelons que l'on a fix\'e un Frob\'enius $Fr$
dans $W_F$. On a $\psi _f(Fr)^k=1$ pour un certain entier $k>0$.
Le centralisateur de $\psi (W_F)$ est donc contenu dans celui du
groupe cyclique engendr\'e par $\psi _{nr}(Fr)^k=\psi (Fr)^k$. Le
nombre de composantes connexes de $^LG$ \'etant fini, il existe un
entier positif $m$, tel que $\psi _{nr}(Fr)^{km}\in\widehat {G}$.

Comme le centralisateur d'un \'el\'ement semi-simple dans un
groupe r\'eductif connexe est r\'eductif, le centralisateur de
$\psi _{nr}(Fr)^{km}$ dans $^LG$ est un sous-groupe r\'eductif
$\Cal {G}$ de $^LG$. Comme $\psi _{nr}(Fr)$ est contenu dans $\Cal
{G}$, l'automorphisme int\'erieur d\'efini par $\psi _{nr}(Fr)$
induit un automorphisme de $\Cal {G}$. Il est d'ordre fini. Comme
$\psi _f(W_F)$ est \'egalement contenu dans $\Cal {G}$, le
centralisateur de $\psi (W_F)$ dans $^LG$ est \'egal au groupe des
points fixes des automorphismes int\'erieurs de $\Cal {G}$
d\'efinis par les \'el\'ements semi-simples $\psi _{nr}(\gamma
)\psi _f(\gamma )$, $\gamma\in W_F$. Or, par ce qui pr\'ec\'edait,
ces automorphismes forment par restriction un sous-groupe fini du
groupe des automorphismes du groupe r\'eductif connexe complexe
$\Cal {G}^{\circ }$. Il est bien connu (cf. par exemple \cite
{PY}) que le groupe des points fixes d'un tel groupe fini
d'automorphismes est un groupe r\'eductif. Ceci prouve la
proposition. \hfill {\fin 2}

\null {\bf 4.2} {\bf Proposition:} \it Soit $\psi : W_F\times \SL
_2(\Bbb C)\rightarrow\ ^LG$ un homomorphisme admissible et soit
$^LM_{\psi }$ un sous-groupe de Levi minimal de $^LG$ qui contient
$\psi (W_F)$. Le tore maximal $\widehat{T}_{\psi }$ contenu dans
le centre de $^LM_{\psi }$ est un tore maximal du centralisateur
de $\psi (W_F)$ dans $^LG$.

\null Preuve: \rm  Le tore maximal $\widehat {T}_{\psi }$ contenu
dans le centre de $^LM_{\psi }$ est inclus dans le centralisateur
de $\psi (W_F)$. Il reste \`a voir qu'il est maximal. Notons
$\widehat {T}^{\psi }$ un tore maximal du centralisateur de $\psi
(W_F)$ qui le contient. Le centralisateur de $\widehat {T}^{\psi
}$ dans $^LG$ contient $\psi (W_F)$. Il se projette donc sur $\Gal
(K/F)$. Par suite, c'est un sous-groupe de Levi de $^LG$ (cf.
\cite {B, 3.5}) qui est n\'ecessairement contenu dans $^LM_{\psi
}$. Par minimalit\'e de $^LM_{\psi }$, on a l'\'egalit\'e, d'o\`u
$\widehat {T}_{\psi }=\widehat {T}^{\psi }$.\hfill {\fin 2}

\null {\bf 4.3} Nous fixons pour la suite de cette section un
sous-groupe parabolique standard $P=MU$ de $G$.

Nous dirons qu'une repr\'esentation irr\'eductible cuspidale
unitaire $(\sigma ,E)$ de $M$ v\'erifie l'hypoth\`ese $(LM)$, si
l'assertion suivante est v\'erifi\'ee.

\null \it (LM) On peut associer \`a $\sigma $ un homomorphisme
admissible discret $\psi _{\sigma }:W_F\times\SL_2(\Bbb
C)\rightarrow\ ^LM$ ayant les propri\'et\'es suivantes: tout
d'abord on suppose que $\psi _{\sigma }$ ait \'et\'e choisi tel
que l'on puisse trouver un sous-groupe de Levi minimal
$^LM_{\sigma }$ de $^LM$ qui contienne $\psi _{\sigma }(W_F)$ et
$s_{\sigma }:=\psi _{\sigma }(\pmatrix q^{1/2} & 0\\ 0 &
q^{-1/2}\endpmatrix )$ et qui soit semi-standard. Notons
$^LM^{\sigma }$ le centralisateur de $\psi _{\sigma }(W_F)$ dans
$^LG$ et $\wt {M}^{\sigma }$ sa composante neutre. Fixons ensuite
une forme int\'erieure $G_0$ de $G$ tel que $^LM_{\sigma }$ soit
un sous-groupe de Levi admissible pour $G_0$. Notons $M_0$ le
sous-groupe de Levi semi-standard de $G_0$ qui correspond \`a
$^LM$ et $M_{\sigma }$ celui qui correspond \`a $^LM_{\sigma }$.
Fixons un sous-groupe parabolique $P_{\sigma }$ de Levi $M_{\sigma
}$ de $G_0$. Notons $\chi _{\lambda _{\sigma }}$ l'\'el\'ement de
$\X (M_{\sigma })$ qui correspond \`a $s_{\sigma }:=\psi _{\sigma
}(\pmatrix q^{1/2} & 0\\ 0 & q^{-1/2}\endpmatrix )$ par la
correspondance de Langlands pour les tores. (Remarquons que
$\lambda _{\sigma }\in a_M^*$, puisque la partie compacte de
$s_{\sigma }$ est n\'ecessairement triviale.) Alors on demande de
plus que l'on peut choisir $G_0$ tel que les conditions suivantes
soient v\'erifi\'ees:

\null a) il existe une repr\'esentation irr\'eductible cuspidale
unitaire $\sigma _0$ de $M_{\sigma }$, tel que $i_{P_{\sigma }\cap
M_0}^{M_0}(\sigma _0\otimes\chi _{\lambda _{\sigma }})$ poss\`ede
un sous-quotient de carr\'e int\'egrable $\pi _0$, de sorte que,
identifiant $\X (M_0)$ et $\X (M)$ par la correspondance de
Langlands pour les tores, $\pi _0$ et $\sigma $ ont m\^eme
fonction $\mu $;

\null b) pour toute racine r\'eduite $\alpha $ dans $\Sigma
(P_{\sigma })$, on a la propri\'et\'e suivante: remarquons que
$(\wt {M}^{\sigma }\cap\wt{M_{\alpha }})^{\circ }$ est \'egal au
centralisateur connexe de $\psi _{\sigma }(W_F)$ dans $^LM_{\alpha
}$. Pour que $\lambda\mapsto\mu ^{M_{\sigma ,\alpha }}(\sigma
_0\otimes\chi _{\lambda\alpha })$ ait un p\^ole en $\lambda >0$,
il faut et il suffit que $\alpha (q)^{\lambda }$ soit un
\'el\'ement $q$-distingu\'e de $(\wt {M}^{\sigma
}\cap\wt{M_{\alpha }})^{\circ }$ et que ce groupe ne soit pas un
tore.

\null \it Remarque: \rm (i) Si la restriction de $\psi _{\sigma }$
\`a $\SL _2({\Bbb C})$ est triviale, alors $^LM_{\sigma }=\ ^LM$.
On peut donc choisir $G_0=G$, $\sigma _0=\sigma $ et la condition
a) est trivialement v\'erifi\'ee.

Il reste la condition b): celle-ci est motiv\'ee par le fait
qu'une repr\'esentation r\'eductible de $M_{\alpha }$
paraboliquement induite par une repr\'esentation cuspidale non
unitaire de $M$ poss\`ede un unique sous-quotient de carr\'e
int\'egrable.

En effet, si $i_{P\cap M_{\alpha }}^{M_{\alpha
}}(\sigma\otimes\chi _{\lambda\alpha })$ est r\'eductible,
$\lambda >0$, cette repr\'esentaion poss\`ede un unique
sous-quotient de carr\'e int\'egrable, et, notant $\psi _{\sigma
}$ le param\`etre de Langlands de $\sigma $, il doit alors exister
un homomorphisme admissible discret $\psi :W_F\times \SL _2(\Bbb
C)\rightarrow\ ^LM_{\alpha }$ tel que $\psi _{\vert W_F}={\psi
_{\sigma }}_{\vert W_F}$ et que $\psi (\pmatrix q^{1/2} & 0\\
0 & q^{-1/2}\endpmatrix )=\alpha (q)^{\lambda }$ (comparer avec
{\bf 5.3}). Ceci implique la propri\'et\'e b).

Inversement, on d\'eduit d'un homomorphisme admissible discret
$\psi :W_F\times \SL _2(\Bbb C)$ $\rightarrow\ ^LM_{\alpha }$ dont
la restriction \`a $W_F$ est donn\'ee par celle de $\psi _{\sigma
}$ une repr\'esentation de carr\'e int\'egrable de $M_{\alpha }$
qui doit \^etre un sous-quotient de $i_{P\cap M_{\alpha
}}^{M_{\alpha }}(\sigma\otimes\chi _{\lambda\alpha })$ pour un
certain nombre r\'eel $\lambda >0$ (comparer avec {\bf 5.4}).

(ii) Remarquons que, quitte \`a remplacer $\psi _{\sigma }$ par un
homomorphisme admissible qui lui est \'equivalent, on peut
toujours trouver un sous-groupe de Levi minimal de $^LM$ qui
contient $\psi _{\sigma }(W_F)$ et qui est semi-standard. Il est
alors uniquement d\'etermin\'e par $\psi _{\sigma }$, puisque
l'intersection de deux sous-groupes de Levi semi-standard est
contenu dans un sous-groupe de Levi semi-standard propre (de
chacun des deux sous-groupes de Levi). Si $^LM_{\sigma }$ est
admissible pour $G$, on peut choisir $G_0=G$. Si la restriction de
$\psi _{\sigma }$ sur $\SL _2({\Bbb C})$ est non triviale, il est
suppos\'e qu'il existe une repr\'esentation de carr\'e
int\'egrable $\pi _0$ de $M$ correspondant au m\^eme param\`etre
$\psi _{\sigma }$ et que celle-ci peut \^etre obtenue \`a partir
d'une repr\'esentation cuspidale $\sigma _0$ de $M_{\sigma }$
comme d\'ecrit ci-dessus. Comme $\pi _0$ et $\sigma $
correspondent au m\^eme param\`etre $\psi _{\sigma }$, leurs
fonctions $\mu $ sont suppos\'ees \^etre \'egales \cite{S,
paragraphe 9}. La condition b) s'explique alors comme dans (i).

(iii) En g\'en\'eral, on peut toujours choisir pour $G_0$ l'unique
forme int\'erieure de $G$ qui est quasi-d\'eploy\'ee sur $F$: tout
sous-groupe de Levi de $^LG$ est admissible pour $G_0$. Par la
correspondance de Langlands pour les formes int\'erieures (qui
n'est connue dans une certaine mesure que pour $\ul{G}_0$ un
groupe lin\'eaire g\'en\'eral), on associe \`a $\sigma $ une
repr\'esentation de carr\'e int\'egrable $\pi _0$ de $M_0$. Les
fonctions $\mu $ de $\sigma $ et $\pi _0$ devraient \^etre
\'egales \cite{S, paragraphe 9}. En fait, $\pi _0$ correspondrait
au m\^eme param\`etre $\psi _{\sigma }$ (cette fois pris relatif
au groupe $M_0$). Les autres conditions s'expliquent comme dans
(ii) (en partant de $\pi _0$).

(iv) Il est possible que l'on puisse \'etablir l'hypoth\`ese
ci-dessus dans une certaine mesure dans le cas o\`u la
localisation des points de r\'eductibilit\'e est connue (cf.
\cite{S}) ou la d\'eduire d'autres propri\'et\'es conjecturales de
ces points de r\'eductibilit\'e \cite{M}, \cite{Z}.

(v) On pourrait penser \`a renforcer la condition b) dans
l'hypoth\`ese (LM). D\'esig-nons pour un nombre complexe $\lambda
$ par $\wt {M}^{\sigma ,\lambda\alpha }$ le centralisateur connexe
de l'image de $W_F$ par l'application $W_F\rightarrow\ ^LM_{\alpha
}$, $\gamma \mapsto \alpha (q)^{v_F(\gamma )\Im (\lambda )}\psi
_{\sigma }(\gamma )$. Si on remplace la derni\`ere phrase dans la
condition b) par "Pour que $\lambda\mapsto\mu ^{M_{\sigma ,\alpha
}}(\sigma _0\otimes\chi _{\lambda\alpha })$ ait un p\^ole en \it
un nombre complexe \rm $\lambda $, il faut et il suffit que
$\alpha (q)^{\Re(\lambda )}$ soit un \'el\'ement $q$-distingu\'e
de $\wt {M}^{\sigma ,\lambda\alpha }$ et que ce groupe ne soit pas
un tore", alors l'hypoth\`ese (LM) serait en quelque sorte
invariante par torsion par un caract\`ere non ramifi\'e unitaire,
i.e., si $\sigma $ v\'erifie $(LM)$, alors, pour tout $\lambda $
dans $a_M^*$, $\sigma\otimes\chi _{\sqrt{-1}\lambda }$
v\'erifierait \'egalement $(LM)$.

En effet, notant $s_{\lambda }$ l'\'el\'ement de $T_{^LM}$ qui
correspond \`a $\lambda $ par la correspondance de Langlands pour
les tores, on peut prendre comme param\`etre de Langlands de
$\sigma\otimes\chi _{\sqrt{-1}\lambda }$ l'homomorphisme
admissible $\psi _{\sigma ,\sqrt{-1}\lambda }:W_F\times\SL_2(\Bbb
C)\rightarrow\ ^LM$, $(\gamma ,h)\mapsto (s_{\lambda
})^{\sqrt{-1}v_F(\gamma )}\psi (\gamma ,h)$. Il est imm\'ediat
qu'il a la propri\'et\'e demand\'ee. (Il faudrait \'eventuellement
faire l'hypoth\`ese suppl\'ementaire que l'\'egalit\'e $\psi
_{\sigma ,\sqrt{-1}\lambda }=\psi _{\sigma }$ \'equivaut \`a dire
que $\sigma $ et $\sigma\otimes\chi _{\sqrt{-1}\lambda }$ sont
isomorphes. Ceci fait partie des conjectures de Langlands pour les
repr\'esentations cuspidales.)

\null {\bf 4.4} Pour \'etablir des r\'esultats d'unicit\'e nous
aurons besoin de l'hypoth\`ese suivante:

\null \it (LMU) Soient $M_1$ et $M_2$ deux sous-groupes de Levi
semi-standard de $G$. Si $\sigma _1$ et $\sigma _2$ sont des
repr\'esentations irr\'eductibles cuspidales unitaires de $M_1$ et
de $M_2$ respectivement qui v\'erifient $(LM)$ et qui sont
conjugu\'ees par un \'el\'ement de $G$, alors leurs param\`etres
de Langlands sont conjugu\'es (\`a l'int\'erieur de $^LG$) par un
\'el\'ement de $\wt {G}$. \rm

\null \it Remarque: \rm Si on conjugue le param\`etre de Langlands
d'une repr\'esentation irr\'e-ductible cuspidale unitaire $\sigma
$ de $M$ \`a l'int\'erieur de $^LG$ par un \'el\'ement de
$\wt{G}$, on obtient un homomorphisme admissible relatif \`a un
autre sous-groupe de Levi semi-standard $M_1$ de $^LG$. Cet
homomorphisme admissible devrait \^etre le param\`etre de
Langlands d'une repr\'esentation irr\'eductible cuspidale unitaire
$\sigma _1$ de $M_1$. Il r\'esulte alors de l'invariance de la
fonction $\mu $ de Harish-Chandra par conjugaison (cf. \cite{W})
que cette repr\'esentation cuspidale $\sigma _1$ doit \'egalement
v\'erifier l'hypoth\`ese $(LM)$.

\null {\bf 4.5} On suppose l'hypoth\`ese $(LM)$ et on garde les
notations pr\'ec\'edentes. On d\'esignera par $\wt {T}_{\sigma }$
le tore maximal dans le centre de $^LM_{\sigma }$. Notons $\Sigma
_{\sigma _0}$ le sous-ensemble de $\Sigma (A_{M_{\sigma }})$
form\'e des racines r\'eduites $\alpha $ tel que $\mu ^{M_{\sigma
,\alpha }}(\sigma_0)=0$. Rappelons que c'est un syst\`eme de
racines \cite{Si2, 3.5}.

\null {\bf Proposition:} \it Le syst\`eme de racines $\Sigma
^{\sigma }(\wt{T}_{\sigma })$ de $\widehat {M}^{\sigma }$ relatif
\`a $\widehat{T}_{\sigma }$ est isomorphe \`a $\Sigma _{\sigma
_0}^{\vee }$, sauf \'eventuellement si $\Sigma _{\sigma _0}$
poss\`ede des facteurs de type $B_n$ ou $C_n$. Dans ce cas,
$\Sigma ^{\sigma }(\wt{T}_{\sigma })$ se d\'ecompose toujours en
facteurs irr\'eductibles selon la d\'ecomposition de $\Sigma
_{\sigma _0}$, mais les facteurs irr\'eductibles de $\Sigma
^{\sigma }(\wt{T}_{\sigma })$ qui correspondent \`a des facteurs
de type $B_n$ ou $C_n$ de $\Sigma _{\sigma _0}$ peuvent aussi bien
\^etre de type $B_n$ que $C_n$.

\null Preuve: \rm Les racines de $\widehat{T}_{\sigma }$ dans
l'alg\`ebre de Lie de $\widehat {M}^{\sigma }$ forment un
sous-ensemble $\Sigma '$ de $\Sigma (\widehat {T}_{\sigma
})=\Sigma (A_{M_{\sigma }})^{\vee }$. Par la condition b) de
l'hypoth\`ese $(LM)$, une racine $\beta ^{\vee }$ appartient \`a
$\Sigma '$, si et seulement si $\beta $ est le multiple d'une
racine r\'eduite $\alpha $ dans $\Sigma (A_{M_{\sigma }})$, telle
que $\chi\mapsto\mu ^{M_{\sigma ,\alpha }}(\sigma_0\otimes\chi )$
ait un p\^ole $\chi _{\lambda }$ avec $\lambda $ \it r\'eel. \rm
Par les propri\'et\'es de la fonction $\mu $ de Harish-Chandra
(cf. {\bf 0.7}), ceci \'equivaut \`a dire que $\mu ^{M_{\sigma
,\alpha }}(\sigma_0)=0$, i.e. $\alpha\in\Sigma _{\sigma _0}$.

Le syst\`eme de racines $\Sigma ^{\sigma }(\wt{T}_{\sigma })$ est
donc form\'e de multiples de racines dans $\Sigma _{\sigma
_0}^{\vee }$. En particulier, $\Sigma ^{\sigma }(\wt{T}_{\sigma
})$ se d\'ecompose en facteurs irr\'eductibles selon la
d\'ecomposition de $\Sigma _{\sigma _0}$. Par ailleurs, la
classification des syst\`emes de racines irr\'eductibles
r\'eduites \cite{Bo} montre qu'un syst\`eme de racines
irr\'eductibles, form\'e de multiples d'un autre syst\`eme de
racines irr\'eductibles, ne peut en diff\'erer en type que si ce
dernier syst\`eme est de type $B_n$ ou $C_n$. Et alors ce
syst\`eme de racines ne peut \^etre que de type $B_n$ ou $C_n$.
\hfill{\fin 2}

\null {\bf Corollaire:} \it Pour que l'on ait l'\'egalit\'e $C(\wt
{M}^{\sigma })^{\circ }=C(^LG)^{\circ }$, il faut et il suffit que
le rang semi-simple de $\wt{M}^{\sigma }$ soit \'egal au rang
parabolique de $^LM_{\sigma }$. Ce nombre est une borne maximale
pour le rang semi-simple de $\wt{M}^{\sigma }$.

\null Preuve: \rm On a  $C(\wt {M}^{\sigma })^{\circ }=(\bigcap
_{\alpha }\ker\alpha )^{\circ}$, o\`u $\alpha $ parcourt les
racines de $\widehat {M}^{\sigma }$ relatives \`a
$\widehat{T}_{\sigma }$. Par la proposition ci-dessus et les faits
r\'esum\'es dans {\bf 1.4}, ce dernier groupe est \'egal \`a
$T_{^LG}$, si et seulement si le rang de $\Sigma _{\sigma _0}$ est
\'egal \`a $\rg(\wt{T}_{\sigma })-\rg(T_{^LG})$. Ce nombre est la
borne maximale pour le rang de $\Sigma _{\sigma _0}$. Ceci prouve
le corollaire. \hfill{\fin 2}

\null {\bf 5.} Continuons \`a fixer jusqu'\`a la fin de {\bf 5.2}
un sous-groupe parabolique standard $P=MU$ de $G$. En outre fixons
une repr\'esentation irr\'eductible cuspidale unitaire $(\sigma
,E)$ de $M$ v\'erifiant l'hypoth\`ese $(LM)$. On notera $\psi
_{\sigma }: W_F \times \SL_2(\Bbb C)\rightarrow\ ^LM$
l'homomorphisme admissible associ\'e \`a $(\sigma ,E)$, et on
garde les notations du paragraphe pr\'ec\'edent relatives \`a
$\psi _{\sigma }$. En particulier, $\wt{M}^{\sigma }$ d\'esignera
la composante neutre du centralisateur de $\psi _{\sigma }(W_F)$
dans $^LG$, et $^LM_{\sigma }$ le sous-groupe de Levi minimal de
$^LM$ contenant $\psi _{\sigma }(W_F)$ qui est semi-standard. En
outre, pour $\lambda\in a_{M,\Bbb C}^*$, on notera $s_{\lambda }$
l'\'el\'ement semi-simple de $T_{^LM}$ qui lui correspond par la
correspondance de Langlands pour les tores (cf. {\bf 1.8}).

\null {\bf 5.1} La propri\'et\'e de $\psi _{\sigma }$ dont on aura
besoin pour prouver notre th\'eor\`eme principal est r\'esum\'ee
dans le lemme ci-dessous. On la d\'eduira de l'hypoth\`ese $(LM)$.

\null {\bf Lemme:} \it Soit $\alpha $ une racine r\'eduite dans
$\Sigma (P)$, $\lambda\in\Bbb R$ et posons $s_{\sigma ,\lambda
\alpha }=s_{\sigma }\alpha (q)^{\lambda }$. Notons $\frak u
_{\pm\alpha ^{\vee }}$ la somme directe des alg\`ebres de Lie des
sous-groupes unipotents de $\widehat {M}^{\sigma }$ associ\'ees
\`a des poids de $T_{^LM}$ \'egaux \`a un multiple non nul de
$\alpha ^{\vee }$. (Cet espace peut \'eventuellement \^etre
r\'eduit \`a z\'ero.)

Consid\'erons $\Ad (s_{\sigma ,\lambda \alpha })$ comme un
endomorphisme de $\frak u _{\pm\alpha ^{\vee }}$. Alors, la
diff\'erence des dimensions des espaces propres associ\'es aux
valeurs propres $q$ et $1$ de $\Ad (s_{\sigma ,\lambda \alpha })$
(cette dimension \'etant \'egale \`a $0$, si $q$ resp. $1$ n'est
pas une valeur propre) est \'egale \`a l'ordre du p\^ole de $\mu
^{M_{\alpha }}(\sigma\otimes\chi )$ en $\chi =\chi _{\lambda\alpha
}$. Plus pr\'ecis\'ement, notant $\frak u _{\pm\alpha ^{\vee
},\lambda \alpha }(z )$ l'espace propre associ\'e \`a une valeur
propre $z$ (et $\frak u _{\pm\alpha ^{\vee },\lambda\alpha }(z
)=0$, si $z$ n'est pas une valeur propre), on a $$\dim (\frak u
_{\pm\alpha ^{\vee },\lambda\alpha }(q))-\dim (\frak u _{\pm\alpha
^{\vee },\lambda\alpha }(1))=\cases
-2,&\text{si $\mu ^{M_{\alpha }}(\sigma\otimes\chi _{\lambda\alpha})=0$}\\
1, &\text{si $\sigma\otimes\chi _{\lambda\alpha }$ est un
p\^ole de $\mu ^{M_{\alpha }}$}\\ 0, &\text{sinon.}\\
\endcases$$

\null Preuve: \rm Avec les notations de {\bf 4.3}, on d\'eduit de
l'hypoth\`ese $(LM)$ et de la formule du produit pour la fonction
$\mu $ de Harish-Chandra que $$\eqalign {&\mu ^{M_{\alpha
}}(\sigma\otimes\chi_{\lambda\alpha })\cr =&\mu ^{M_{0,\alpha
}}(\pi _0\otimes\chi_{\lambda\alpha })\cr =&{(\mu ^{M_{0,\alpha
}}/\mu ^{M_0})}(\sigma _0\otimes\chi _{\lambda _{\sigma }}\otimes
\chi _{\lambda\alpha })\cr =&\prod _{\beta }\mu ^{M_{\sigma ,\beta
}}(\sigma _0\otimes\chi _{\lambda _{\sigma }}\otimes\chi
_{\lambda\alpha }),\cr }$$ o\`u $\beta $ parcourt les racines
r\'eduites dans $\Sigma (P_{\sigma })$ de restriction \`a $A_M$
\'egale \`a un multiple de $\alpha $.

Notons $\Sigma ^{\sigma }(\wt{T}_{\sigma })$ le syst\`eme de
racines de $\wt{M}^{\sigma }$ relatif au tore maximal
$\wt{T}_{\sigma }$ dans le centre de $\wt{M}_{\sigma }$. La
condition b) de l'hypoth\`ese $(LM)$ et les propositions {\bf 4.2}
et {\bf 4.5} (y inclus la preuve de {\bf 4.5}) impliquent que $\mu
^{M_{\sigma ,\beta }}(\sigma _0\otimes\chi _{\lambda _{\sigma
}}\otimes\chi )$ a un p\^ole en $\chi =\chi _{\lambda\alpha }$, si
et seulement si $\beta \in\Sigma _{\sigma_0}$ et si toute racine
$\gamma $ dans $\Sigma ^{\sigma }(\wt{T}_{\sigma })$ qui est un
multiple de $\beta ^{\vee }$ v\'erifie $\gamma (s_{\sigma ,\lambda
})\in\{q^{-1},q\}$. Les p\^oles de $\mu $ relatifs \`a une
repr\'esentation cuspidale sont n\'ecessairement d'ordre 1 (cf.
{\bf 0.7}). Il r\'esulte par ailleurs des propri\'et\'es de la
fonction $\mu $ que $\mu ^{M_{\sigma ,\beta }}(\sigma _0
\otimes\chi _{\lambda _{\sigma }}\otimes\chi _{\lambda\alpha
})=0$, si et seulement si $\Re (\langle\beta ^{\vee },\lambda
_{\sigma }+\lambda\alpha\rangle )=0$ (ce qui \'equivaut \`a $\vert
\beta ^{\vee }(s_{\sigma ,\lambda\alpha })\vert =1$), et si $\mu
^{M_{\sigma ,\beta }}(\sigma _0 \otimes \chi _{\lambda _{\sigma
}}\otimes\chi )$ a un p\^ole en $\chi =\chi _{\lambda\alpha
+\lambda '\beta }$ pour un certain r\'eel $\lambda '>0$. Ceci
\'equivaut par ce qui pr\'ec\'edait  \`a dire que $\beta\in\Sigma
_{\sigma_0}$ et que $\gamma (s_{\sigma ,\lambda })=1$ pour toute
racine $\gamma $ dans $\Sigma ^{\sigma }(\wt{T}_{\sigma })$ qui
est un multiple de $\beta ^{\vee }$. Les z\'eros de la fonction
$\mu $ relatifs \`a une repr\'esentation cuspidale sont
n\'ecessairement d'ordre $2$ (cf. {\bf 0.7}). \hfill {\fin 2}

\null {\bf 5.2 Proposition:} \it Soit $\lambda \in a_M^{G*}$. Pour
que $\sigma\otimes\chi _{\lambda }$ soit un p\^ole de la fonction
$\mu $ de Harish-Chandra d'ordre \'egal au rang parabolique de
$M$, il faut et il suffit que le rang semi-simple de
$\wt{M}^{\sigma }$ soit \'egal au rang parabolique de $^LM_{\sigma
}$ et que $s_{\sigma ,\lambda }=s_{\sigma }s_{\lambda }$ soit
$q$-distingu\'e dans $\wt{M}^{\sigma }$.

\null Preuve: \rm   Remarquons d'abord que le centralisateur
connexe de $\psi _{\sigma }(W_F)$ dans $^LM$ est \'egal \`a
$(\wt{M}\cap\wt{M}^{\sigma })^{\circ }$. Par {\bf 4.2}, le tore
maximal $\widehat {T}_{\sigma }$ contenu dans le centre de
$^LM_{\sigma }$ est un tore maximal de $\widehat {M}^{\sigma }$.

Notons $\frak u_{\widehat {M}^{\sigma }, \pm}$ (resp. $\frak
u_{\widehat {M}^{\sigma }\cap \widehat {M}, \pm}$) le sous-espace
vectoriel de l'alg\`ebre de Lie de $\wt{M}^{\sigma }$ (resp.
$(\widehat {M}^{\sigma }\cap \widehat {M})^{\circ }$) engendr\'e
par les \'el\'ements nilpotents. Consid\'erons $Ad (s_{\sigma
,\lambda })$ comme endomorphisme de ces espaces vectoriels. Notons
$\frak u_.(q)$ et $\frak u_.(1)$ respectivement les espaces
propres associ\'es aux valeurs propres $q$ et $1$ (ces espaces
\'etant \'egaux \`a $0$, si $q$ resp. $1$ n'est pas une valeur
propre). Remarquons que $s_{\sigma ,\lambda }$ est $q$-distingu\'e
dans $\wt{M}^{\sigma }$, si et seulement si $\dim (\frak
u_{\widehat {M}^{\sigma }, \pm}(q))-\dim (\frak u_{\widehat
{M}^{\sigma }, \pm}(1))=\rg _{ss}(\widehat {M}^{\sigma })$.

Tenant compte de l'in\'egalit\'e $\rg _{ss}(\wt{M}^{\sigma })\leq
\rg_{ss}(^LG)-\rg _{ss}(^LM_{\sigma })$ (cf. {\bf 4.5}) et de la
remarque dans {\bf 3.5}, tout revient \`a montrer que
$$\dim (\frak u_{\widehat {M}^{\sigma }, \pm}(q))-\dim (\frak
u_{\widehat {M}^{\sigma }, \pm}(1))=\rg _{ss}(^LG)-\rg
_{ss}(^LM_{\sigma })\tag 5-2-1$$ \'equivaut \`a
$$\ord _{\sigma\otimes\chi _{\lambda }}\mu =\rg _{ss}(G)-\rg
_{ss}(M).\tag 5-2-2$$

Avec les notations de {\bf 5.1}, on a $$\frak u_{\widehat
{M}^{\sigma }, \pm}=\frak u_{\widehat {M}^{\sigma }\cap \widehat
{M}, \pm}\oplus \bigoplus_{\alpha ^{\vee }} \frak u_{\pm\alpha
},\tag 5-2-3$$ o\`u $\alpha ^{\vee }$ parcourt les racines
r\'eduites dans $\Sigma (P)$.

On d\'eduit du lemme {\bf 5.1} et de la formule du produit pour la
fonction $\mu $ que
$$\sum _{\alpha ^{\vee }}(\dim (\frak u_{\pm\alpha}(q))-\dim (\frak u_{\pm\alpha}(1)))=
\ord _{\sigma\otimes\chi _{\lambda }}\mu .\tag 5-2-4$$

Supposons d'abord $\widehat {M_{\sigma }}=\widehat {M}$, ce qui
implique que le rang parabolique de $^LM_{\sigma }$ est \'egal \`a
celui de $M$ (cf. {\bf 1.4}). On a $\frak u_{\widehat {M}^{\sigma
}\cap \widehat {M}, \pm}=0$, et on d\'eduit de (5-2-3) et (5-2-4)
que
$$\dim (\frak u_{\widehat {M}^{\sigma }, \pm}(q))-\dim (\frak
u_{\widehat {M}^{\sigma }, \pm}(1))=\ord _{\sigma\otimes\chi
_{\lambda }}\mu .$$ L'\'equivalence de (5-2-1) et de (5-2-2) est
alors imm\'ediate.

Supposons maintenant $\widehat {M_{\sigma }}\ne\widehat {M}$. Par
l'hypoth\`ese $(LM)$, il existe une forme int\'erieure $G_0$ de
$G$, des sous-groupes de Levi semi-standard $M_{\sigma }$ et $M_0$
de $G_0$, $M_{\sigma }\subseteq M_0$, tels que $M_0$ et $M_{\sigma
}$ correspondent respectivement aux sous-groupes de Levi $^LM$ et
$^LM_{\sigma }$ de $^LG$, ainsi qu'une repr\'esentation cuspidale
$\sigma _0$ de $M_0$ et $\lambda _{\sigma }\in a_{M_{\sigma
}}^{M_0*}$, tels que, identifiant $\X (M_0)$ et $\X (M)$ par la
correspondance de Langlands pour les tores, on ait $(\mu /\mu
^{M_0})(\sigma\otimes\chi _{\lambda _{\sigma }}\otimes\chi )=\mu
(\sigma \otimes\chi )$ pour tout $\chi\in \X (M)$. Par ailleurs,
fixant un sous-groupe parabolique $P_{\sigma }$ de $G_0$ de
facteur de Levi $M_{\sigma }$, $i_{P_{\sigma }\cap
M_0}^{M_0}(\sigma _0\otimes\chi _{\lambda _{\sigma }})$ poss\`ede
un sous-quotient de carr\'e int\'egrable et donc, par le
r\'esultat principal de \cite{H2} rappel\'e en {\bf 0.7}, $\mu
^{M_0}$ a un p\^ole d'ordre $\rg _{ss}(M_0)-\rg_{ss}(M_{\sigma })$
en $\sigma _0\otimes\chi _{\lambda _{\sigma }}$. De plus, le
param\`etre de Langlands de $\sigma _0$ est donn\'e par
l'homomorphisme admissible $W_F\rightarrow\ ^LM_{\sigma }$,
$\gamma\mapsto\psi _{\sigma }(\gamma ,1)$, et il v\'erifie
l'hypoth\`ese $(LM)$.

Le cas trait\'e ci-dessus, appliqu\'e \`a $M_0$, $\sigma _0$ et
$\lambda _{\sigma }$, nous montre que $$\dim (\frak u_{\widehat
{M}^{\sigma }\cap M, \pm}(q))-\dim (\frak u_{\widehat {M}^{\sigma
}\cap M, \pm}(1))=\rg _{ss}(^LM)-\rg _{ss}(^LM_{\sigma }).\tag
5-2-5$$ En addition (5-2-5) et (5-2-4) et en utilisant (5-2-3), on
trouve que le c\^ot\'e gauche de (5-2-1) vaut $\rg _{ss}(^LM)-\rg
_{ss}(^LM_{\sigma })+\ord _{\sigma \otimes\chi _{\lambda }}\mu $.
L'\'equivalence de (5-2-1) et de (5-2-2) est alors une
cons\'equence de la proposition dans {\bf 1.4} et de l'\'egalit\'e
\'el\'ementaire
$$\rg_{ss}(^LG)-\rg_{ss}(^LM_{\sigma })=(\rg_{ss}(^LG)-\rg_{ss}(^LM))
+(\rg_{ss}(^LM)-\rg_{ss}(^LM_{\sigma }))$$ \hfill {\fin 2}

\null {\bf 5.3} Soit $\pi $ une repr\'esentation de carr\'e
int\'egrable de $G$. Il existe alors un sous-groupe parabolique
semi-standard $P=MU$ de $G$, une repr\'esentation irr\'eductible
cuspidale unitaire $\sigma $ de $M$ et $\lambda\in a_M^*$,
$\lambda _G=0$, tel que $\pi $ soit un sous-quotient de
$i_P^G(\sigma\otimes\chi _{\lambda })$. Par le r\'esultat
principal de \cite{H2, 8.2} rappel\'e en {\bf 0.7},
$\sigma\otimes\chi _{\lambda }$ est un p\^ole de la fonction $\mu
$ de Harish-Chandra d'ordre \'egal au rang parabolique de $M$.

Supposons que $\sigma $ v\'erifie l'hypoth\`ese $(LM)$. On sait
donc lui associer un para-m\`etre de Langlands $\psi _{\sigma
}:W_F\times\SL _2(\Bbb C)\rightarrow\ ^LM$ avec les propri\'et\'es
\'enonc\'ees dans $(LM)$. Alors, par {\bf 5.2}, le rang
semi-simple de $\wt{M}^{\sigma }$ est \'egal au rang parabolique
de $^LM_{\sigma }$, et $s_{\sigma ,\lambda }:=s_{\sigma
}s_{\lambda }$ est un \'el\'ement $q$-distingu\'e du
centralisateur connexe $\wt {M}^{\sigma }$ de $\psi _{\sigma
}(W_F)$ dans $^LG$. Il doit donc \^etre contenu dans le groupe
d\'eriv\'e de $\wt {M}^{\sigma }$. La partie compacte dans la
d\'ecomposition polaire de $s_{\sigma ,\lambda }$ est triviale,
puisque il en est ainsi pour $s_{\sigma }$ et $s_{\lambda }$.
Choisissons \`a l'aide de {\bf 3.6} un \'el\'ement nilpotent
$N_{\sigma ,\lambda }$ de l'alg\`ebre de Lie de $\wt{M}^{\sigma }$
tel que $(s_{\sigma ,\lambda },N_{\sigma ,\lambda })$ soit une
$L^2$-paire dans $\wt{M}^{\sigma }$.

Comme remarqu\'e dans {\bf 3.7}, on peut en d\'eduire un
homomorphisme admissible et discret $\phi _{\sigma ,\lambda }:\SL
_2(\Bbb C)\rightarrow\ \wt{M}^{\sigma }$, v\'erifiant $\phi
_{\sigma ,\lambda }($ $\pmatrix 1 & 1 \\ 0 & 1\endpmatrix )=\exp
(N_{\lambda, \sigma })$ et $\phi _{\sigma ,\lambda }(\pmatrix
q^{1/2}& 0\\ 0 &q^{-1/2}\\\endpmatrix)=s_{\sigma ,\lambda }$.

On pose $$\psi _{\pi }(\gamma ,h)=\psi _{\sigma }(\gamma ,1)\phi
_{\sigma ,\lambda }(h).$$

\null {\bf Th\'eor\`eme:} \it L'application $\psi _{\pi
}:W_F\times \SL_2(\Bbb C)\rightarrow\ ^LG$, $(\gamma
,h)\mapsto\psi _{\pi }(\gamma ,h)$, ci-dessus est bien d\'efinie.
C'est un homomorphisme admissible et discret. Si on admet de plus
l'hypoth\`ese $(LMU)$ de {\bf 4.4}, alors $\psi _{\pi }$ est, \`a
conjugaison pr\`es, uniquement d\'etermin\'e par $\pi $.

\null Preuve: \rm L'application $\psi _{\pi }$ est \'evidemment
bien d\'efinie. Montrons que son image n'est contenue dans aucun
sous-groupe de Levi propre de $^LG$. Le centralisateur de l'image
de $\psi _{\pi }$ est \'egal \`a l'intersection du centralisateur
de $\psi _{\sigma }(W_F)$ et de celui de l'image de $\phi _{\sigma
,\lambda }$. Par cons\'equence, tout tore de $^LG$ qui centralise
l'image de $\psi _{\pi }$ doit \^etre contenu dans $\wt{M}^{\sigma
}$ et centraliser l'image de $\phi _{\sigma ,\lambda }$. Or, comme
remarqu\'e ci-dessus, $\phi _{\sigma ,\lambda }$ est discret. Un
tel tore doit donc \^etre inclus dans le centre de $\wt
{M}^{\sigma }$. Mais, par le corollaire de {\bf 4.5}, comme le
rang semi-simple de $\wt{M}^{\sigma }$ est \'egal au rang
parabolique de $^LM_{\sigma }$, ceci implique qu'il est inclus
dans le centre de $^LG$. Il ne peut donc y avoir de sous-groupe de
Levi propre de $^LG$ qui contient l'image de $\psi _{\pi }$.

Ceci prouve que $\psi _{\pi}$ est discret, et, en particulier, que
la propri\'et\'e (v) dans la d\'efinition {\bf 1.5} d'un
homomorphisme admissible est v\'erifi\'ee. Les autres
propri\'et\'es de cette d\'efinition sont imm\'ediates.

Quant \`a l'unicit\'e, on sait par la th\'eorie des
repr\'esentations des groupes $p$-adiques que le triplet $(\sigma
,M,\lambda )$ est, \`a conjugaison par un \'el\'ement de $G$
pr\`es, uniquement d\'etermin\'e par $\pi $. Soit $(\sigma ',
M',\lambda ')$ un autre triplet qui lui est conjugu\'e par un
\'el\'ement de $G$. Alors le param\`etre de Langlands $\psi
_{\sigma '}$ de $\sigma '$ est par l'hypoth\`ese $(LMU)$
conjugu\'e \`a $\psi _{\sigma }$ par un \'el\'ement de $\wt {G}$.
On peut modifier cet \'el\'ement convenablement, pour qu'il
conjugue \'egalement $s_{\lambda }$ et $s_{\lambda '}$ ainsi que
$\exp(N_{\sigma ,\lambda })$ et $\exp(N_{\sigma ',\lambda '})$,
puisque ces \'el\'ements sont contenus dans le centralisateur de
$\psi _{\sigma }(W_F)$ et de $\psi _{\sigma '}(W_F)$
respectivement. Il est alors imm\'ediat que l'homomorphisme
admissible discret que l'on d\'eduit de $(\sigma ', M',\lambda ')$
par le proc\'ed\'e ci-dessus est \'equivalent \`a $\psi _{\pi }$.
\hfill{\fin 2}

\null {\bf 5.4} Pour montrer la r\'eciproque de {\bf 5.3}, i.e.
que tout homomorphisme admissible discret provient d'une
repr\'esentation de carr\'e int\'egrable de $G$ comme d\'ecrit
dans {\bf 5.3}, nous avons besoin de l'hypoth\`ese
suppl\'ementaire suivante qui caract\'erise les homomorphismes
admissibles correspondant \`a des $L$-paquets form\'es de
repr\'esentations cuspidales. Cette hypoth\`ese sera justifi\'ee
post\'erieurement (voir la remarque ci-dessous).

\null \it (LC) Soit $M$ un sous-groupe de Levi semi-standard de
$G$ et soit $\psi :W_F\times\SL _2(\Bbb C)\rightarrow\ ^LM$ un
homomorphisme admissible discret. Supposons que $^LM$ soit un
sous-groupe de Levi admissible minimal de $^LG$ tel que $\psi
(W_F)$ soit contenu dans $^LM$ et que $\psi (1,\pmatrix q^{1/2} & 0\\
0 & q^{-1/2}\endpmatrix)$ soit un \'el\'ement $q$-distingu\'e du
centralisateur connexe de $\psi (W_F)$ dans $^LM$.

Alors, il existe une repr\'esentation irr\'eductible cuspidale
unitaire $\sigma $ de $M$ telle que $\psi $ soit le param\`etre de
Langlands de $\sigma $ et que $\sigma $ et $\psi $ v\'erifient
$(LM)$ et $(LMU)$. \rm

\null Nous avons alors

\null {\bf Th\'eor\`eme:} \it Supposons l'hypoth\`ese $(LC)$
v\'erifi\'ee pour tout sous-groupe de Levi semi-standard de $G$.
Alors, tout homomorphisme admissible discret $\psi :W_F\times\SL
_2(\Bbb C)\rightarrow\ ^LG$ est \'equivalent \`a un homomorphisme
admissible $\psi _{\pi }$ pour une certaine repr\'esentation de
carr\'e int\'egrable $\pi $ de $G$.

\null Preuve: \rm Comme $\psi $ est discret, $T_{^LG}$ doit \^etre
le tore maximal contenu dans le centre du centralisateur connexe
$\wt{M}^{\psi }$ de $\psi (W_F)$. On d\'eduit alors de {\bf 4.2}
que le rang semi-simple de $\wt{M}^{\psi }$ est \'egal au rang
parabolique d'un sous-groupe de Levi minimal $^LM_{\psi }$ de
$^LG$ qui contient $\psi (W_F)$. Choisissons un sous-groupe de
Levi admissible minimal $^LM$ de $^LG$ qui contient $\psi (W_F)$
et qui est tel que $s:=\psi (1,\pmatrix q^{1/2} & 0\\0 &
q^{-1/2}\endpmatrix)$ soit un \'el\'ement $q$-distingu\'e du
centralisateur connexe de $\psi (W_F)$ dans $^LM$. (Si $^LM_{\psi
}$ est admissible, alors on peut \'evidemment poser $^LM=\
^LM_{\psi }$, tout \'el\'ement d'un tore complexe \'etant
trivialement $q$-distingu\'e.)

Quitte \`a conjuguer $\psi $, on peut choisir $^LM$ tel que
$^LM\supseteq\ ^LM_{\psi }\supseteq\ ^LT$ et que $s$ soit contenu
dans le tore maximal $\wt{T}_{\psi }$ contenu dans le centre de
$^LM_{\psi }$. Notons $M$ le sous-groupe de Levi semi-standard de
$G$ qui correspond \`a $^LM$. On a d\'ej\`a remarqu\'e que la
partie compacte de $s$ dans la d\'ecomposition polaire doit \^etre
triviale. On peut donc \'ecrire $s=s_Ms^M$ avec $s_M\in T_{^LM}$
et $s^M\in \wt {T}_{\psi }\cap\wt{M}^{der}$ tel que $\chi
(s_M)=\chi (s)$ pour tout caract\`ere alg\'ebrique
$^LM\rightarrow{\Bbb C}^{\times }$. Par choix de $M$,
l'\'el\'ement $s^M$ est $q$-distingu\'e dans $(\wt {M}^{\psi
}\cap\ ^LM)^{\circ }$. Le rang semi-simple de $(\wt {M}^{\psi
}\cap\ ^LM)^{\circ }$ est \'egal au rang parabolique de
$^LM_{\sigma }$ relatif \`a $^LM$, puisque $^LM$ est semi-standard
et que le rang semi-simple de $\wt{M}^{\psi }$ est \'egal au rang
parabolique de $^LM_{\psi }$ relatif \`a $^LG$. On conclut que
$s^M$ est contenu dans le groupe d\'eriv\'e de $\wt{M}^{\psi }$.

Comme dans {\bf 5.3} (sans d\'etour par un \'el\'ement $\lambda $
de l'alg\`ebre de Lie de $T_{^LM}$), on en d\'eduit un morphisme
admissible discret $W_F\times\SL _2(\Bbb C)\rightarrow\ ^LM$ qui,
par l'hypoth\`ese $(LC)$, est le param\`etre de Langlands d'une
repr\'esentation irr\'eductible cuspidale unitaire $\sigma $ de
$M$. Notons ce morphisme admissible $\psi _{\sigma }$. En
particulier, la restriction de $\psi _{\sigma }$ \`a $W_F$ est
donn\'ee par celle de $\psi $, et on a l'\'egalit\'e $\psi_{\sigma
}(1,\pmatrix q^{1/2} & 0\\0 & q^{-1/2}\endpmatrix)=s^M$. Notons
$\lambda $ l'\'el\'ement de $a_M^*$ qui correspond \`a $s_M$ par
la correspondance de Langlands pour les tores.

Comme $\sigma $ et $\psi _{\sigma }$ v\'erifient par hypoth\`ese
l'hypoth\`ese $(LM)$ et que $s$ est $q$-distingu\'e dans
$\wt{M}^{\psi }$, on d\'eduit de la proposition {\bf 5.2} que
$\sigma\otimes\chi _{\lambda }$ est un p\^ole de la fonction $\mu
$ de Harish-Chandra d'ordre \'egal au rang parabolique de $M$. Il
suit alors du r\'esultat principal de \cite{H2} rappel\'e en {\bf
0.7} que $i_P^G(\sigma\otimes\chi _{\lambda })$ poss\`ede un
sous-quotient de carr\'e int\'egrable $\pi $. Quitte \`a choisir
un \'el\'ement nilpotent convenable dans la construction de $\psi
_{\pi }$ dans {\bf 5.3}, on obtient bien $\psi =\psi _{\pi }$.
\hfill{\fin 2}

\null \it Remarque: \rm (i) Le th\'eor\`eme justifie
post\'erieurement l'hypoth\`ese $(LC)$: en effet, les
homomorphismes admissibles discrets qui y sont caract\'eris\'es
sont ceux qui ne peuvent pas correspondre \`a une repr\'esentation
de carr\'e int\'egrable obtenue par induction parabolique. Ces
param\`etres doivent donc correspondre \`a des repr\'esenta-tions
cuspidales.

Remarquons qu'il existe bien des cas o\`u $\psi (1,\pmatrix q^{1/2} & 0\\
0 & q^{-1/2}\endpmatrix)$ n'est pas $q$-distin-gu\'e relatif au
sous-groupe de Levi admissible minimal qui contient $\psi (W_F)$
(par ex. $\psi (W_F)=\{1\}$ et $G$ un certain groupe non
quasi-d\'eploy\'e de type $B_n$).

\null (ii) Si on veut de plus que le $L$-paquet de $\pi $ soit
uniquement d\'etermin\'e par $\psi $, il faut ajouter
l'hypoth\`ese suivante:

\null \it (LCU) (a) Soit $\psi :W_F\times\SL _2(\Bbb
C)\rightarrow\ ^LM$ un homomorphisme admissible discret. Alors
deux sous-groupes de Levi admissibles $^LM_1$ et $^LM_2$ de $^LG$
qui sont minimaux pour la propri\'et\'e qu'ils contiennent $\psi
(W_F)$ et que $s_{\psi }:=\psi (1,\pmatrix q^{1/2} & 0\\
0 & q^{-1/2}\endpmatrix)$ soit $q$-distingu\'e dans le
centralisateur connexe de $\psi (W_F)$ dans $^LM_1$ (resp.
$^LM_2$) sont conjugu\'es par un \'el\'ement de $\wt{G}$ qui
centralise $\psi (W_F)$ et $s_{\psi }$.

(b) Supposons que $M$ et $M'$ soient deux sous-groupes de Levi
semi-standard de $G$ et que $\psi :W_F\times\SL _2(\Bbb
C)\rightarrow\ ^LM$, $\psi ':W_F\times\SL _2(\Bbb C)\rightarrow\
^LM'$ soient deux homomorphismes admissibles qui sont conjugu\'es
par un \'el\'ement de $\wt {G}$. Alors, si $\psi $ est le
param\'etre de Langlands d'une repr\'esentation irr\'eductible
cuspidale unitaire $\sigma $ de $M$ et que les propri\'et\'es
$(LM)$ et $(LMU)$ sont v\'erifi\'ees, alors $\psi '$ est le
param\'etre de Langlands d'une repr\'esentation irr\'eductible
cuspidale unitaire $\sigma '$ de $M'$ et $\sigma $ et $\sigma '$
sont conjugu\'es par un \'el\'ement de $G$. En particulier, les
propri\'et\'es $(LM)$ et $(LMU)$ sont v\'erifi\'ees relatives \`a
$\sigma '$. \rm

\null La partie a) de l'hypoth\`ese (LCU) devrait se laisser
v\'erifier \`a la main (il faut comparer les diagrammes des
paraboliques distingu\'es avec ceux des syst\`emes de racines des
groupe r\'eductifs non quasi-d\'eploy\'es), mais nous avons
effectu\'e cette v\'erification seulement dans des cas
particuliers, repoussant le cas g\'en\'eral \`a un autre moment.

\null {\bf 5.5} Admettons toutes nos hypoth\`eses ($(LM)$,
$(LMU)$, $(LC)$ et $(LCU)$(i),(ii)) et terminons en d\'ecrivant
sommairement un proc\'ed\'e qui associe \`a tout homomorphisme
admissible $W_F\times\SL _2(\Bbb C)\rightarrow\ ^LG$ (le
$L$-paquet d') une repr\'esentation irr\'eductible lisse de $G$,
ainsi qu'un proc\'ed\'e inverse associant \`a toute
repr\'esentation irr\'eductible lisse de $G$ un homomorphisme
admissible $W_F\times\SL _2(\Bbb C)\rightarrow\ ^LG$.

Notons $^LM'$ un sous-groupe de Levi minimal de $^LG$ contenant
l'image de $\psi $. Il est uniquement d\'etermin\'e \`a
conjugaison pr\`es par un \'el\'ement du centralisateur de $\psi
(W_F\times\SL _2(\Bbb C))$ (cf. \cite{B, 3.6}). Quitte \`a
remplacer $\psi $ par un homomorphisme admissible qui lui est
\'equivalent, on peut supposer $^LM'$ semi-standard.
L'homomorphisme $\psi ^{M'}:W_F \times \SL _2(\Bbb C)\rightarrow\
^LM'$, obtenu par restriction \`a droite de $\psi $ se d\'ecompose
en un produit $\psi _{M',nr}\psi _{rc}^{M'}$, o\`u $\psi
_{rc}^{M'}:W_F \times \SL _2(\Bbb C)\rightarrow\ ^LM'$ est discret
et o\`u $\psi _{M',nr}:W_F \times \SL _2(\Bbb C)\rightarrow\ ^LM'$
est trivial sur $\SL _2(\Bbb C)$, non ramifi\'e sur $W_F$ et \`a
image dans le groupe des \'el\'ements semi-simples hyperboliques
de la composante neutre du centre de $^LM'$ (ceci r\'esulte de
{\bf 4.1}, en utilisant la d\'ecomposition polaire d'un
\'el\'ement semi-simple).

Par {\bf 5.4} et $(LMU)$, $\psi _{rc}^{M'}$ est \'egal \`a un
homomorphisme $\psi _{\tau _0}^{M'}$ associ\'e \`a une
repr\'esentation irr\'eductible de carr\'e int\'egrable $\tau _0$
de $M'$.

Notons $\lambda _{M'}$ un \'el\'ement de $a_{M'}^*$ tel que $\chi
_{\lambda _{M'}}$ corresponde \`a $\psi _{M',nr}(Fr)$ par la
correspondance de Langlands pour les tores et $M''$ le sous-groupe
de Levi semi-standard de $G$ qui contient $M'$ et qui est
engendr\'e par les racines $\alpha $ dans $\Sigma (A_{M'})$ qui
v\'erifient $\langle\alpha^{\vee },\Re (\lambda _{M'})\rangle =0$.
Notons $\lambda _{M'}=\lambda _{M'}^{M''}\oplus\lambda _{M''}$ la
d\'ecomposition de $\lambda _{M'}$ selon la d\'ecomposition
$a_{M'}^*=a_{M'}^{M''*}\oplus a_{M''}^*$. Il existe un sous-groupe
parablique $P''$ de $G$ de facteur de Levi $M''$ tel que $\lambda
_{M''}>_{P''}0$. Choisissons un sous-groupe parabolique $P'$ de
$M''$ de facteur de Levi $M'$. La repr\'esentation induite
$i_{P'}^{M''}\tau _0$ est une somme directe de repr\'esentations
temp\'er\'ees. Choisissons-en une sous-repr\'esentation
irr\'eductible $\tau $.

On d\'efinit alors $\pi _{\psi }$ comme \'etant la
repr\'esentation irr\'eductible lisse d\'eduite de $P'',\tau,
\lambda _{M''}$ par la classification de Langlands. Le $L$-paquet
de $\pi _{\psi }$ est par construction uniquement d\'etermin\'e
par $\psi $.

\null D\'ecrivons le proc\'ed\'e inverse de celui d\'ecrit
ci-dessus, en admettant toutes nos hypoth\`eses: soit $\pi $ une
repr\'esentation irr\'eductible lisse de $G$. Par la
classification de Langlands, celle-ci est associ\'ee \`a un
sous-groupe parabolique semi-standard $P''=M''U''$ de $G$, une
repr\'esentation temp\'er\'ee $\tau $ de $M''$ et un \'el\'ement
$\lambda _{M''}$ de $a_{M''}^*$ qui est strictement positif dans
la chambre de Weyl associ\'ee \`a $P''$. Ces donn\'ees sont
uniquement d\'etermin\'ees \`a conjugaison pr\`es. On peut trouver
un sous-groupe parabolique semi-standard $P'=M'U'$ de $M''$ et une
repr\'esentation de carr\'e int\'egrable $\tau _0$ de $M'$ tels
que $\tau $ soit une sous-repr\'esentation de $i_{P'}^{M''}\tau
_0$. Par {\bf 5.3}, on d\'eduit de $\tau _0$ un homomorphisme
admissible discret $\psi _{\tau _0}^{M'}:W_F\times\SL _2(\Bbb
C)\rightarrow\ ^LM'$. Notons $s_{M''}$ l'\'el\'ement semi-simple
dans $C(^LM'')^{\circ }$ qui correspond \`a $\lambda _{M''}$ par
la correspondance de Langlands pour les tores. Alors on d\'efinit
$\psi _{\pi }:W_F\times\SL _2(\Bbb C)\rightarrow\ ^LG$ par
$(\gamma ,h)\mapsto s_{M''}^{v_F(\gamma )}\psi _{\tau
_0}^{M'}(\gamma ,h)$.

\null Le r\'esultat suivant est une cons\'equence imm\'ediate des
constructions:

\null {\bf Th\'eor\`eme:} \it Pour tout homomorphisme admissible
$\psi :W_F\times \SL _2(\Bbb C)\rightarrow\ ^LG$, $\psi _{\pi
_{\psi }}$ est \'equivalent \`a $\psi $. Inversement, pour toute
repr\'esentation irr\'eductible lisse $\pi $ de $G$, $\pi _{\psi
_{\pi }}$ appartient au m\^eme $L$-paquet que $\pi $. \hfill \fin
2

\enddocument
\bye